\documentclass[12pt]{article}
\usepackage{amsmath}

\usepackage{t1enc}
\usepackage[latin1]{inputenc}
\usepackage[english]{babel}
\usepackage{amsmath,amsthm}
\usepackage{amsfonts}
\usepackage{latexsym}
\usepackage[dvips]{graphicx}
\usepackage{graphicx}
\usepackage[natural]{xcolor}
\newtheorem{theorem}{Theorem}
\newtheorem{lemma}[theorem]{Lemma}
\newtheorem{proposition}[theorem]{Proposition}

\newtheorem{example}[theorem]{Example}
\newtheorem{definition}[theorem]{Definition}

\begin{document}

\title{{\LARGE Maximum Principle for Forward-Backward Doubly Stochastic Control
Systems and Applications }}

\author{Liangquan Zhang$^{1,}$$^2$ and Yufeng SHI$^1$
\thanks{This work was partially Supported by National Natural Science Foundation of China Grant 10771122,
Natural Science Foundation of Shandong Province of China Grant
Y2006A08 and National Basic Research Program of China (973 Program,
No. 2007CB814900)}
\thanks {E-mail: yfshi@sdu.edu.cn (Y. Shi), Liangquan.Zhang@etudiant.univ-brest.fr(L. Zhang).}\\
{\footnotesize {(1.School of Mathematics, Shandong University} }\\
{\footnotesize {2. Laboratoire de Math\'ematiques Universit\'e de
Bretagne
Occidentale)} }\\
{\footnotesize {}}}
 \maketitle

\date{}
\maketitle

\begin{abstract}
\baselineskip 0.65cm

The maximum principle for optimal control problems of fully coupled
forward-backward doubly stochastic differential equations (FBDSDEs
in short) in the global form is obtained, under the assumptions that
the diffusion coefficients do not contain the control variable, but
the control domain need not to be convex. We apply our stochastic
maximum principle (SMP in short) to investigate the optimal control
problems of a class of stochastic partial differential equations
(SPDEs in short). And as an example of the SMP, we solve a kind of
forward-backward doubly stochastic linear quadratic optimal control
problems as well. In the last section, we use the solution of
FBDSDEs to get the explicit form of the optimal control for linear
quadratic stochastic optimal control problem and open-loop Nash
equilibrium point for nonzero sum differential games problem.



\end{abstract}

\section{Introduction}
\hspace{0.25in}It is well known that optimal control problem is one
of the central themes of control science. The necessary conditions
of optimal problem were established for deterministic control system
by Pontryakin's group [24] in the 1950's and 1960's. Since then, a
lot of work has been done on the forward stochastic system such as
Kushner 13], Bismut [5], Bensoussan [2, 3], Haussmann [9, 10] and
Peng [20] etc.

Peng [20] studied the following type of stochastic optimal control
problem. Minimize a cost function\'e
\[
J\left( v_{\left( \cdot \right) }\right) ={\bf E}\int_0^Tl\left(
x_t,v_t\right) \text{d}t+{\bf E}\left( h_T\right) ,
\]
subject to
\begin{equation}
\left\{
\begin{array}{l}
\text{d}x_t=g\left( t,x_t,v_t\right) \text{d}t+\sigma \left(
t,x_t,v_t\right) \text{d}B_t, \\
x_0=x,
\end{array}
\right.   \tag{1.1}
\end{equation}
over an admissible control domain which need not be convex, and the
diffusion coefficients contain the control variable. In his paper,
by spike variational method and the second order adjoint equations,
Peng [20] obtained a general stochastic maximum principle for the
above optimal control problem. It was just the adjoint equations in
stochastic optimal control problems that motivated the famous theory
of backward stochastic differential equations (BSDEs in short) (see
[18]). Later Peng [21] studied a stochastic optimal control problem
where state variables are described by the system of forward and
backward SDEs, that is
\begin{equation}
\left\{
\begin{array}{l}
\text{d}x_t=f\left( t,x_t,v_t\right) \text{d}t+\sigma \left(
t,x_t,v_t\right) \text{d}W_t, \\
x_0=x, \\
\text{d}y_t=g\left( t,x_t,v_t\right) \text{d}t+z_t\text{d}W_t, \\
y_T=y,
\end{array}
\right.   \tag{1.2}
\end{equation}
where $x$ and $y$ are given deterministic constants. The optimal
control problem is to minimize the cost function
\[
J\left( v_{\left( \cdot \right) }\right) ={\bf E}\left[
\int_0^Tl\left( t,x_t,y_t,v_t\right) \text{d}t+h\left( x_T\right)
+\gamma \left( y_0\right) \right] ,
\]
over an admissible control domain which is convex. Xu [28] studied
the following non-fully coupled forward-backward stochastic control
system
\begin{equation}
\left\{
\begin{array}{l}
\text{d}x_t=f\left( t,x_t,v_t\right) \text{d}t+\sigma \left(
t,x_t\right)
\text{d}W_t, \\
x_0=x, \\
\text{d}y_t=g\left( t,x_t,y_t,z_t,v_t\right) \text{d}t+z_t\text{d}W_t, \\
y_T=h\left( x_T\right) .
\end{array}
\right.   \tag{1.3}
\end{equation}
The optimal control problem is to minimize the cost function
\[
J\left( v_{\left( \cdot \right) }\right) ={\bf E}\gamma \left(
y_0\right) ,
\]
over ${\cal U}_{ad},$ but the control domain is non-convex. Wu [26]
firstly gave the maximum principle for optimal control problem of
fully coupled forward-backward stochastic system
\begin{equation}
\left\{
\begin{array}{l}
\text{d}x_t=f\left( t,x_t,y_t,z_t,v_t\right) \text{d}t+\sigma \left(
t,x_t,y_t,z_t,v_t\right) \text{d}B_t, \\
\text{d}y_t=-g\left( t,x_t,y_t,z_t,v_t\right) \text{d}t+z_t\text{d}B_t, \\
x_0=x,\quad \quad y_T=\xi ,
\end{array}
\right.   \tag{1.4}
\end{equation}
where $\xi $ is a random variable and the cost function
\[
J\left( v_{\left( \cdot \right) }\right) ={\bf E}\left[
\int_0^TL\left( t,x_t,y_t,z_t,v_t\right) \text{d}t+\Phi \left(
x_T\right) +h\left( y_0\right) \right] .
\]
The optimal control problem is to minimize the cost function
$J\left( v_{\left( \cdot \right) }\right) $ over an admissible
control domain which is convex. Ji and Zhou [12] obtained a maximum
principle for stochastic optimal control of non-fully coupled
forward-backward stochastic system with terminal state constraints.
Shi and Wu [25] studied the maximum principle for fully coupled
forward-backward stochastic system
\begin{equation}
\left\{
\begin{array}{l}
\text{d}x_t=b\left( t,x_t,y_t,z_t,v_t\right) \text{d}t+\sigma \left(
t,x_t,y_t,z_t\right) \text{d}B_t, \\
\text{d}y_t=-f\left( t,x_t,y_t,z_t,v_t\right) \text{d}t+z_t\text{d}B_t, \\
x_0=x,\quad \quad y_T=h\left( x_T\right) .
\end{array}
\right.   \tag{1.5}
\end{equation}
and the cost function is
\[
J\left( v_{\left( \cdot \right) }\right) ={\bf E}\left[
\int_0^Tl\left( t,x_t,y_t,z_t,v_t\right) \text{d}t+\Phi \left(
x_T\right) +\gamma \left( y_0\right) \right] .
\]
The control domain is non-convex but the forward diffusion does not
contain the control variable. For more details in this field, see
Yong and Zhou [29].

In order to provide a probabilistic interpretation for the solutions
of a class of semilinear stochastic partial differential equations
(SPDEs in short), Pardoux and Peng [19] introduced the following
backward doubly stochastic differential equation (BDSDE in short):
\begin{eqnarray*}
Y_t=\xi +\int\limits_t^Tf(s,Y_s,Z_s)\text{d}s+\int\limits_t^Tg(s,Y_s,Z_s)%
\text{\^d}B_s-\int\limits_t^TZ_s\text{d}W_s,\quad 0\leq t\leq T.
\end{eqnarray*}
\begin{equation}
\tag{1.6}
\end{equation}
Note that the integral with respect to $\{B_t\}$ is a ``backward
It\^o integral'' and the integral with respect to $\{W_t\}$ is a
standard forward It\^o integral. These two types of integrals are
particular cases of the It\^o-Skorohod integral (for details see
[15]). Peng and Shi [22] introduced a type of time-symmetric
forward-backward stochastic differential equations, i.e., so-called
fully coupled forward-backward doubly stochastic differential
equations (FBDSDEs in short):
\begin{equation}
\left\{
\begin{array}{lll}
y_t & = & x+\int\limits_0^tf\left( s,y_s,Y_s,z_s,Z_s\right) \text{d}%
s+\int\limits_0^tg\left( s,y_s,Y_s,z_s,Z_s\right) \text{d}%
W_s-\int\limits_0^tz_s\text{\^d}B_s, \\
Y_t & = & h\left( y_T\right) +\int\limits_t^TF\left(
s,y_s,Y_s,z_s,Z_s\right) \text{d}s+\int\limits_t^TG\left(
s,y_s,Y_s,z_s,Z_s\right)
\text{\^d}B_s+\int\limits_t^TZ_s\text{d}W_s.
\end{array}
\right.   \tag{1.7}
\end{equation}
In FBDSDEs (1.7), the forward equation is ``forward'' with respect
to a standard stochastic integral d$W_t$, as well as ``backward''
with respect to a backward stochastic integral \^d$B_t$; the coupled
``backward equation'' is ``forward'' under the backward stochastic
integral \^d$B_t$ and ``backward'' under the forward one. In other
words, both the forward equation and the backward one are types of
BDSDE (1.6) with different directions of stochastic integrals. So
(1.7) provides a very general framework of fully coupled
forward-backward stochastic systems. Peng and Shi [22] proved the
existence and uniqueness of solutions to FBDSDEs (1.7) with
arbitrarily fixed time duration under some monotone assumptions.
FBDSDEs (1.7) can provide a probabilistic interpretation for the
solutions of a class of quasilinear SPDEs.

As we have known, stochastic control problem of the SPDEs arising
from partial observation control has been studied by Mortensen [9],
using a dynamic programming approach, and subsequently by
Bensoussan, using a maximum principle method. See [4], [16] and the
references therein for more information. Our approach differs from
the one of Bensoussan. More precisely, we relate the FBDSDEs to one
kind of SPDEs with control variables where the control systems of
SPDEs can be transformed to the relevant control systems of FBDSDEs.
To our knowledge, this is the first time to treat the optimal
control problems of SPDEs from a new perspective of FBDSDEs. It is
worth mentioning that the quasilinear SPDEs in [17] \O ksendal
considered can just be related to our partially coupled FBDSDEs.

Besides, in Section 6 we investigate the nonzero sum stochastic
differential game problem. This problem have been considered by
Friedman [8], Bensoussan [1] and Eisele [7]. For stochastic case
Hammadene [11] and Wu [27] (for more information see references
therein) showed existence result of Nash equilibrium point under
some assumptions, respectively. Here, we extend their result to
doubly stochastic case in which we can regard the backward
filtration as the disturbed information come from outside the
''control system''.

In this paper, we consider the following fully coupled
forward-backward doubly stochastic control system
\begin{equation}
\left\{
\begin{array}{lll}
y_t & = & x+\int_0^tf\left( s,y_s,Y_s,z_s,Z_s,v_s\right) \text{d}%
s+\int_0^tg\left( s,y_s,Y_s,z_s,Z_s\right) \text{d}W_s-\int_0^tz_s\text{\^d}%
B_s, \\
Y_t & = & h\left( y_T\right) +\int_t^TF\left(
s,y_s,Y_s,z_s,Z_s,v_s\right)
\text{d}s+\int_t^TG\left( s,y_s,Y_s,z_s,Z_s\right) \text{\^d}B_s+\int_t^TZ_s%
\text{d}W_s.
\end{array}
\right.  \tag{1.8}
\end{equation}
Our optimal control problem is to minimize the cost function:
\[
J\left( v_{\left( \cdot \right) }\right) ={\bf E}\left[
\int_0^Tl\left( t,y_t,Y_t,z_t,Z_t,v_t\right) \text{d}t+\Phi \left(
y_T\right) +\gamma \left( Y_0\right) \right]
\]
over an admissible control domain which need not be convex. It is
obvious that (1.8) covers (1.3) and (1.5), so (1.8) can describe
more intricate control systems. As for the fully coupled
forward-backward doubly stochastic control systems such as (1.8)
whose diffusion coefficients contain the control variables, this
issue will be carried out in our future publications.

The notable difficulties to obtain the maximum principles for the
fully coupled forward-backward doubly stochastic control systems
within non-convex control domains are how to use the spike
variational method to get variational equations with enough high
order estimates and how to use the duality technique to obtain the
adjoint equations. On account of the quadruple of variables in the
FBDSDEs, we can not directly apply the methods used in [25], [26]
and [28]. In this paper, by virtue of the results of FBDSDEs in
[22], we can ensure the existence and uniqueness of the solutions
for the adjoint FBDSDEs which are obtained by applying the duality
technique to the variational equations. Besides, we apply the
technique of FBDSDEs to get the enough high order estimates for the
solutions of the variational equations.

This paper is organized as follows. In Section 2, we state the
problems and some assumptions. In Section 3, we study the
variational equations and variational inequalities. In Section 4, a
stochastic maximum principle in global form is obtained,
subsequently, an example of this kind of control problems is given
in this section. As an application, we study the optimal control
problem of a kind of SPDEs with control variable by the approach of
FBDSDEs in Section 5. Lastly, we give the explicit form of Nash
equilibrium point for a kind of stochastic differential game
problem.

For the simplicity of notations, we only consider the case where
both $y$ and $Y$ are one-dimensional, and the control $v$ is also
one-dimensional. While in order to give the general results, we
consider the multi-dimensional case in Section 6.

\section{Statement of the problem}

\hspace{0.25in}Let $\left( \Omega ,\mathcal{F},P\right) $ be a complete
probability space, and $[0,T]$ be a given time duration throughout this
paper. Let $\left\{ W_t;0\leq t\leq T\right\} $ and $\left\{ B_t;0\leq t\leq
T\right\} $ be two mutually independent standard Brownian motions defined on
$\left( \Omega ,\mathcal{F},P\right) $, with values respectively in $\mathbf{%
R}^d$ and in $\mathbf{R}^l$. Let $\mathcal{N}$ denote the class of $P$-null
elements of $\mathcal{F}$. For each $t\in \left[ 0,T\right] $, we define
\[
\mathcal{F}_t\doteq \mathcal{F}_t^W\vee \mathcal{F}_{t,T}^B
\]
where $\mathcal{F}_t^W=\mathcal{N}\vee \sigma \left\{ W_r-W_0;0\leq r\leq
t\right\} $, $\mathcal{F}_{t,T}^B=\mathcal{N}\vee \sigma \left\{
B_r-B_t;t\leq r\leq T\right\} $. Note that the collection $\left\{ \mathcal{F%
}_t,t\in \left[ 0,T\right] \right\} $ is neither increasing nor decreasing,
and it does not constitute a classical filtration. We introduce the following

\begin{definition}
A stochastic process $X=\left\{ X_t;t\geq 0\right\} $ is called $\mathcal{F}%
_t$-progressively measurable, if for any $t\geq 0$, $X$ on $\Omega \times
\left[ 0,t\right] $ is measurable with respect to $\left( \mathcal{F}%
_t^W\times \mathcal{B}\left( \left[ 0,t\right] \right) \right) \vee \left(
\mathcal{F}_{t,T}^B\times \mathcal{B}\left( \left[ t,T\right] \right)
\right) $.
\end{definition}

Let $M^2\left( 0,T;\mathbf{R}^n\right) $ denote the space of all (classes of
$dP\otimes dt$ a.e. equal) $\mathbf{R}^n$-valued $\mathcal{F}_t$%
-progressively measurable stochastic processes $\left\{ v_t;t\in \left[
0,T\right] \right\} $ which satisfy
\[
\mathbf{E}\int_0^T\left| v_t\right| ^2\text{d}t<\infty .
\]
Obviously $M^2\left( 0,T;\mathbf{R}^n\right) $ is a Hilbert space. For a
given $u\in M^2\left( 0,T;\mathbf{R}^d\right) $ and $v\in M^2\left( 0,T;%
\mathbf{R}^l\right) $, one can define the (standard) forward It\^o's
integral $\int_0^{\cdot }u_s$d$W_s$ and the backward It\^o's integral $%
\int_{\cdot }^Tv_s$\^d$B_s$. They are both in $M^2\left( 0,T;\mathbf{R}%
\right) $, (see [14] for details).

Let $L^2\left( \Omega ,\mathcal{F}_T,P;\mathbf{R}\right) $ denote the space
of all $\mathcal{F}_T$-measurable one-valued random variable $\xi $
satisfying $\mathbf{E}\left| \xi \right| ^2<\infty .$ Under this framework,
we consider the following forward-backward doubly stochastic control system.
\begin{equation}
\left\{
\begin{array}{l}
dy_t=f\left( t,y_t,Y_t,z_t,Z_t,v_t\right) \text{d}t+g\left(
t,y_t,Y_t,z_t,Z_t\right) \text{d}W_t-z_t\text{\^d}B_t, \\
dY_t=-F\left( t,y_t,Y_t,z_t,Z_t,v_t\right) \text{d}t-G\left(
t,y_t,Y_t,z_t,Z_t\right) \text{\^d}B_t+Z_t\text{d}W_t, \\
y_0=x,\quad \quad Y_T=h\left( y_T\right) ,\quad \quad t\in \left[ 0,T\right]
,
\end{array}
\right.  \tag{2.1}
\end{equation}
where $\left( y_{\left( \cdot \right) },Y_{\left( \cdot \right) },z_{\left(
\cdot \right) },Z_{\left( \cdot \right) },v_{\left( \cdot \right) }\right)
\in \mathbf{R\times R\times R}^l\mathbf{\times R}^d\times \mathbf{R},$ $x\in
\mathbf{R}$ is a given constant$,$ $T>0,$%
\begin{eqnarray*}
F &:&\left[ 0,T\right] \times \mathbf{R}\times \mathbf{R}\times \mathbf{R}%
^l\times \mathbf{R}^d\times \mathbf{R}\rightarrow \mathbf{R,} \\
f &:&\left[ 0,T\right] \times \mathbf{R}\times \mathbf{R}\times \mathbf{R}%
^l\times \mathbf{R}^d\times \mathbf{R}\rightarrow \mathbf{R,} \\
G &:&\left[ 0,T\right] \times \mathbf{R}\times \mathbf{R}\times \mathbf{R}%
^l\times \mathbf{R}^d\times \mathbf{R}\rightarrow \mathbf{R}^l, \\
g &:&\left[ 0,T\right] \times \mathbf{R}\times \mathbf{R}\times \mathbf{R}%
^l\times \mathbf{R}^d\times \mathbf{R}\rightarrow \mathbf{R}^d, \\
h &:&\mathbf{R}\rightarrow \mathbf{R.}
\end{eqnarray*}
Let $\mathcal{U}$ be a nonempty subset of $\mathbf{R.}$ We define the
admissible control set
\[
\mathcal{U}_{ad}\doteq \left\{ v_{\left( \cdot \right) }\in M^2\left( 0,T;%
\mathbf{R}\right) ;\text{ }v_t\in \mathcal{U},\text{ }0\leq t\leq T,\text{
a.e., a.s.}\right\} .
\]
Our optimal control problem is to minimize the cost function:
\begin{equation}
J\left( v_{\left( \cdot \right) }\right) \doteq \mathbf{E}\left[
\int_0^Tl\left( t,y_t,Y_t,z_t,Z_t,v_t\right) \text{d}t+\Phi \left(
y_T\right) +\gamma \left( Y_0\right) \right]  \tag{2.2}
\end{equation}
over $\mathcal{U}_{ad}$, where
\begin{eqnarray*}
l &:&\left[ 0,T\right] \times \mathbf{R}\times \mathbf{R}\times \mathbf{R}%
^l\times \mathbf{R}^d\times \mathbf{R}\rightarrow \mathbf{R,} \\
\Phi &:&\mathbf{R\rightarrow R,} \\
\gamma &:&\mathbf{R\rightarrow R.}
\end{eqnarray*}
An admissible control $u_{\left( \cdot \right) }$ is called an optimal
control if it attains the minimum over $\mathcal{U}_{ad}$. That is to say,
we want to find a $u_{\left( \cdot \right) }$ such that
\[
J\left( u_{\left( \cdot \right) }\right) \doteq \inf\limits_{v_{\left( \cdot
\right) }\in \mathcal{U}_{ad}}J\left( v_{\left( \cdot \right) }\right) .
\]
(2.1) is called the state equation, the solution $\left(
y_t,Y_t,z_t,Z_t\right) $ corresponding to $u_{\left( \cdot \right) }$ is
called the optimal trajectory.

Next we will give some notations:
\[
\zeta =\left(
\begin{array}{c}
y \\
Y \\
z \\
Z
\end{array}
\right) ,\quad A\left( t,\zeta \right) =\left(
\begin{array}{c}
-F \\
f \\
-G \\
g
\end{array}
\right) \left( t,\zeta \right) .
\]
We use the usual inner product $\left\langle \cdot ,\cdot \right\rangle $
and Euclidean norm $\left| \cdot \right| $ in $\mathbf{R}$, $\mathbf{R}^l$
and $\mathbf{R}^d.$ All the equalities and inequalities mentioned in this
paper are in the sense of $dt\otimes dP$ almost surely on $\left[ 0,T\right]
\times \Omega .$ We assume that

\begin{enumerate}
\item[(H1)]  For each $\zeta \in \mathbf{R}^{1+1+l+d},$ $A\left( \cdot
,\zeta \right) $ is an $\mathcal{F}_t$-measurable process defined on $\left[
0,T\right] $ with $A\left( \cdot ,0\right) \in M^2\left( 0,T;\mathbf{R}%
^{1+1+l+d}\right) .$

\item[(H2)]  $A\left( t,\zeta \right) $ and $h\left( y\right) $ satisfy
Lipschitz conditions: there exists a constant $k>0,$ such that
\[
\left\{
\begin{array}{l}
\left| A\left( t,\zeta \right) -A\left( t,\bar \zeta \right) \right| \leq
k\left| \zeta -\bar \zeta \right| ,\text{\quad }\forall \zeta ,\text{ }\bar
\zeta \in \mathbf{R}^{1+1+l+d},\text{ }\forall t\in \left[ 0,T\right] , \\
\left| h\left( y\right) -h\left( \bar y\right) \right| \leq k\left| y-\bar
y\right| ,\quad \forall y,\text{ }\bar y\in \mathbf{R.}
\end{array}
\right.
\]
\end{enumerate}

The following monotonic conditions introduced in [17], are the main
assumptions in this paper.

\begin{enumerate}
\item[(H3)]  $\left\{
\begin{array}{l}
\left\langle A\left( t,\zeta \right) -A\left( t,\bar \zeta \right) ,\zeta
-\bar \zeta \right\rangle \leq -\mu \left| \zeta -\bar \zeta \right| ^2, \\
\quad \forall \zeta =\left( y,Y,z,Z\right) ,\text{ }\bar \zeta =\left( \bar
y,\bar Y,\bar z,\bar Z\right) \in \mathbf{R\times R\times R}^l\mathbf{\times
R}^d,\text{ }\forall t\in \left[ 0,T\right] . \\
\left\langle h\left( y\right) -h\left( \bar y\right) ,y-\bar y\right\rangle
\geq 0,\text{ }\forall y,\text{ }\bar y\in \mathbf{R,}
\end{array}
\right. $
\end{enumerate}

\noindent or

\begin{enumerate}
\item[(H3)']  $\left\{
\begin{array}{l}
\left\langle A\left( t,\zeta \right) -A\left( t,\bar \zeta \right) ,\zeta
-\bar \zeta \right\rangle \geq \mu \left| \zeta -\bar \zeta \right| ^2, \\
\quad \forall \zeta =\left( y,Y,z,Z\right) ,\text{ }\bar \zeta =\left( \bar
y,\bar Y,\bar z,\bar Z\right) \in \mathbf{R\times R\times R}^l\mathbf{\times
R}^d,\text{ }\forall t\in \left[ 0,T\right] . \\
\left\langle h\left( y\right) -h\left( \bar y\right) ,y-\bar y\right\rangle
\leq 0,\text{ }\forall y,\text{ }\bar y\in \mathbf{R,}
\end{array}
\right. $
\end{enumerate}

\noindent where $\mu $ is a positive constant.

\begin{proposition}
For any given admissible control $v_{\left( \cdot \right) },$ we assume
(H1), (H2) and (H3) (or (H1), (H2) and (H3)') hold. Then FBDSDE (2.1) has a
unique solution $\left( y_t,Y_t,z_t,Z_t\right) \in M^2\left( 0,T;\mathbf{R}%
^{1+1+l+d}\right) .$
\end{proposition}

\noindent The proof is referred to [17]. We need a farther assumption as
follows:

\begin{enumerate}
\item[(H4)]  $F,$ $f,$ $G,$ $g,$ $h,$ $l,$ $\Phi ,$ $\gamma $ are
continuously differentiable with respect to $\left( y,Y,z,Z\right) ,$ $y$
and $Y$. They and all their derivatives are bounded by a constant $C$.
\end{enumerate}

\noindent Lastly, we need the following extension of It\^o's formula (for
details see [14]).

\begin{proposition}
Let
\[
\alpha \in S^2\left( 0,T;\mathbf{R}^k\right) ,\beta \in M^2\left( 0,T;%
\mathbf{R}^k\right) ,\gamma \in M^2\left( 0,T;\mathbf{R}^{k\times l}\right)
,\delta \in M^2\left( 0,T;\mathbf{R}^{k\times d}\right)
\]
satisfy:
\[
\alpha _t=\alpha _0+\int_0^t\beta _s\text{d}s+\int_0^t\gamma _s\text{\^d}%
B_s+\int_0^t\delta _s\text{d}W_s,\quad \text{ }0\leq t\leq T.
\]
Then
\begin{eqnarray*}
\left| \alpha _t\right| ^2=\left| \alpha _0\right| ^2+2\int_0^t\left( \alpha
_s,\beta _s\right) \text{d}s+2\int_0^t\left( \alpha _s,\gamma _s\text{\^d}%
B_s\right) +2\int_0^t\left( \alpha _s,\delta _s\text{d}W_s\right) \\
-\int_0^t\left| \gamma _s\right| ^2\text{d}s+\int_0^t\left| \delta _s\right|
^2\text{d}s, \\
\mathbf{E}\left| \alpha _t\right| ^2=\mathbf{E}\left| \alpha _0\right| ^2+2%
\mathbf{E}\int_0^t\left( \alpha _s,\beta _s\right) \text{d}s-\mathbf{E}%
\int_0^t\left| \gamma _s\right| ^2\text{d}s+\mathbf{E}\int_0^t\left| \delta
_s\right| ^2\text{d}s.
\end{eqnarray*}
\end{proposition}

\noindent Here $S^2\left( 0,T;\mathbf{R}^k\right) $ denotes the space of
(classes of $dP\otimes dt$ a.e. equal) all $\mathcal{F}_t$-progressively
measurable $k$-dimensional processes $v$ with
\[
\mathbf{E}\left( \sup\limits_{0\leq t\leq T}\left| v_t\right| ^2\right)
<\infty .
\]

\section{Variational equations and variational inequalities}

\hspace{0.25in}Suppose $\left( y_t,Y_t,z_t,Z_t,u_t\right) $ is the solution
to our optimal control problem. We introduce the following spike variational
control:
\[
u_t^\varepsilon =\left\{
\begin{array}{l}
v,\quad \tau \leq t\leq \tau +\varepsilon , \\
u_t,\quad \text{otherwise,}
\end{array}
\right.
\]
where $\varepsilon >0$ is sufficiently small, $\tau \in \left[ 0,T\right] $.
$v$ is an arbitrary $\mathcal{F}_\tau $-measurable random variable with
values in $\mathcal{U},$ $0\leq t\leq T,$ and $\sup\limits_{\omega \in
\Omega }\left| v\left( \omega \right) \right| <\infty $ Let $\left(
y_t^\varepsilon ,Y_t^\varepsilon ,z_t^\varepsilon ,Z_t^\varepsilon \right) $
be the trajectory of the control system (2.1) corresponding to the control $%
u_t^\varepsilon .$

For convenience, we use the following notations in this paper:
\begin{eqnarray*}
\Xi _y &=&\Xi _y\left( t,y_t,Y_t,z_t,Z_t,u_t\right) , \\
\Xi _y\left( u_t^\varepsilon \right) &=&\Xi _y\left(
t,y_t,Y_t,z_t,Z_t,u_t^\varepsilon \right) , \\
\Xi \left( u_t\right) &=&\Xi \left( t,y_t,Y_t,z_t,Z_t,u_t\right) , \\
\Xi \left( u_t^\varepsilon \right) &=&\Xi \left(
t,y_t,Y_t,z_t,Z_t,u_t^\varepsilon \right) , \\
&&\text{etc,}
\end{eqnarray*}
where $\Xi =f,$ $F,$ $g,$ $G$, respectively. We introduce the following
variational equations:
\begin{equation}
\left\{
\begin{array}{lll}
\text{d}y_t^1 & = & \left[ f_yy_t^1+f_YY_t^1+f_zz_t^1+f_ZY_t^1+f\left(
u_t^\varepsilon \right) -f\left( u_t\right) \right] \text{d}t \\
&  & +\left[ g_yy_t^1+g_YY_t^1+g_zz_t^1+g_ZZ_t^1\right] \text{d}W_t-z_t^1%
\text{\^d}B_t, \\
y_0^1 & = & 0, \\
\text{d}Y_t^1 & = & -\left[ F_yy_t^1+F_YY_t^1+F_zz_t^1+F_ZZ_t^1+F\left(
u_t^\varepsilon \right) -F\left( u_t\right) \right] \text{d}t \\
&  & -\left[ G_yy_t^1+G_YY_t^1+G_zz_t^1+G_ZZ_t^1\right] \text{\^d}B_t+Z_t^1%
\text{d}W_t, \\
Y_T^1 & = & h_y\left( y_T\right) y_T^1.
\end{array}
\right.  \tag{3.1}
\end{equation}
Owing to (H4), it is easy to check that the variational equation (3.1) same
as (2.1), also satisfies (H1), (H2) and (H3). Thus by Proposition 2, there
exists a unique solution $\left( y_t^1,Y_t^1,z_t^1,Z_t^1\right) \in \mathbf{%
R\times R\times R}^l\mathbf{\times R}^d,$ $0\leq t\leq T,$ satisfying (3.1).
The variational inequalities can be derived from the fact $J\left( u_{\left(
\cdot \right) }^\varepsilon \right) -J\left( u_{\left( \cdot \right)
}\right) \geq 0.$ The following lemmas play important roles to establish the
inequalities.

\begin{lemma}
We assume (H1)-(H4) hold. Then we have
\begin{equation}
\mathbf{E}\int_0^T\left| y_t^1\right| ^2\text{d}t\leq C\varepsilon ,
\tag{3.2}
\end{equation}
\begin{equation}
\mathbf{E}\int_0^T\left| Y_t^1\right| ^2\text{d}t\leq C\varepsilon ,
\tag{3.3}
\end{equation}
\begin{equation}
\mathbf{E}\int_0^T\left| z_t^1\right| ^2\text{d}t\leq C\varepsilon ,
\tag{3.4}
\end{equation}
\begin{equation}
\mathbf{E}\int_0^T\left| Z_t^1\right| ^2\text{d}t\leq C\varepsilon .
\tag{3.5}
\end{equation}

\noindent where $C>0$ is some constant.
\end{lemma}

\noindent \textbf{Proof.}\quad\noindent Using the It\^o's formula to $%
\left\langle y_t^1,Y_t^1\right\rangle ,$ it follows that
\begin{eqnarray*}
&&\mathbf{E}\left( \left| y_T^1\right| ^2h_y\left( y_T\right) \right) \\%
[0.03in]
&=&\mathbf{E}\int_0^T\left( f_yy_t^1+f_YY_t^1+f_zz_t^1+f_ZZ_t^1\right) Y_t^1%
\text{d}t \\
&&-\mathbf{E}\int_0^T\left( F_yy_t^1+F_YY_t^1+F_zz_t^1+F_ZZ_t^1\right) y_t^1%
\text{d}t \\
&&-\mathbf{E}\int_0^T\left( G_yy_t^1+G_YY_t^1+G_zz_t^1+G_ZZ_t^1\right) z_t^1%
\text{d}t \\
&&+\mathbf{E}\int_0^T\left( g_yy_t^1+g_YY_t^1+g_zz_t^1+g_ZZ_t^1\right) Z_t^1%
\text{d}t \\
&&+\mathbf{E}\int_0^T\left( f\left( u_t^\varepsilon \right) -f\left(
u_t\right) \right) Y_t^1\text{d}t \\
&&-\mathbf{E}\int_0^T\left( F\left( u_t^\varepsilon \right) -F\left(
u_t\right) \right) y_t^1\text{d}t.  \label{aas}
\end{eqnarray*}
\begin{equation}
\tag{3.6}
\end{equation}

Since (3.1) satisfies the monotonic condition (H3), it is easy to see that
\begin{eqnarray*}
&&\mathbf{E}\left( \left| y_T^1\right| ^2h_y\left( y_T\right) \right) +\mu
\mathbf{E}\int_0^T\left( \left| y_t^1\right| ^2+\left| Y_t^1\right|
^2+\left| z_t^1\right| ^2+\left| Z_t^1\right| ^2\right) \text{d}t  \nonumber
\\
&\leq &\mathbf{E}\int_0^T\left( f\left( u_t^\varepsilon \right) -f\left(
u_t\right) \right) Y_t^1\text{d}t-\mathbf{E}\int_0^T\left( F\left(
u_t^\varepsilon \right) -F\left( u_t\right) \right) y_t^1\text{d}t  \nonumber
\\
&\leq &\frac 1\mu \mathbf{E}\int_0^T\left| f\left( u_t^\varepsilon \right)
-f\left( u_t\right) \right| ^2\text{d}t+\frac \mu 4\mathbf{E}\int_0^T\left|
Y_t^1\right| ^2\text{d}t  \nonumber \\
&&+\frac 1\mu \mathbf{E}\int_0^T\left| F\left( u_t^\varepsilon \right)
-F\left( u_t\right) \right| ^2\text{d}t+\frac \mu 4\mathbf{E}\int_0^T\left|
y_t^1\right| ^2\text{d}t.
\end{eqnarray*}
\begin{equation}
\tag{3.7}
\end{equation}

From (H4) and (3.7), it is easy to know that (3.2)-(3.5) hold. The proof is
complete. \quad $\Box$

However, the order of the estimate for $\left(
y_t^1,Y_t^1,z_t^1,Z_t^1\right) $ is too low to get the variational
inequalities. We need to give some more elaborate estimates. For that, we
firstly give the following lemma.

\begin{lemma}
Assuming (H1)-(H4) hold, then we have
\begin{equation}
\sup\limits_{0\leq t\leq T}\left( \mathbf{E}\left| y_t^1\right| ^2\right)
\leq C\varepsilon ,  \tag{3.8}
\end{equation}
\begin{equation}
\sup\limits_{0\leq t\leq T}\left( \mathbf{E}\left| Y_t^1\right| ^2\right)
\leq C\varepsilon .  \tag{3.9}
\end{equation}
\end{lemma}

\noindent \textbf{Proof.}\quad \noindent Squaring both sides of
\begin{eqnarray*}
y_t^1+\int_0^tz_s^1\text{\^d}B_s &=&\int_0^t\left(
f_yy_s^1+f_YY_s^1+f_zz_s^1+f_ZZ_s^1+f\left( u_s^\varepsilon \right) -f\left(
u_s\right) \right) \text{d}s \\
&&+\int_0^t\left( g_yy_s^1+g_YY_s^1+g_zz_s^1+g_ZZ_s^1\right) \text{d}W_s,
\end{eqnarray*}
noting that
\[
\mathbf{E}\left[ y_t^1\int_0^tz_s^1\text{\^d}B_s\right] =\mathbf{E}\left[
\mathbf{E}^{\mathcal{F}_t}\left( y_t^1\int_0^tz_s^1\text{\^d}B_s\right)
\right] =\mathbf{E}\left[ y_t^1\mathbf{E}^{\mathcal{F}_t}\left( \int_0^tz_s^1%
\text{\^d}B_s\right) \right] =0,
\]
we have
\begin{eqnarray*}
&&\mathbf{E}\left| y_t^1\right| ^2+\mathbf{E}\int_0^t\left| z_s^1\right| ^2%
\text{d}s \\
&=&\mathbf{E[}\int_0^t\left( f_yy_s^1+f_YY_s^1+f_zz_s^1+f_ZZ_s^1+f\left(
u_s^\varepsilon \right) -f\left( u_s\right) \right) \text{d}s \\
&&+\int_0^t\left( g_yy_s^1+g_YY_s^1+g_zz_s^1+g_ZZ_s^1\right) \text{d}W_s]^2
\\
&\leq &C\mathbf{E}\int_0^t\left[ \left| y_s^1\right| ^2+\left| Y_s^1\right|
^2+\left| z_s^1\right| ^2+\left| Z_s^1\right| ^2\right] \text{d}s \\
&&+C\mathbf{E}\left( \int_0^t\left( f\left( u_s^\varepsilon \right) -f\left(
u_s\right) \right) \text{d}s\right) ^2.
\end{eqnarray*}
Thus
\[
\sup\limits_{0\leq t\leq T}\left( \mathbf{E}\left| y_t^1\right| ^2\right)
\leq C\varepsilon .
\]

By the similar argument, we can have
\[
\sup\limits_{0\leq t\leq T}\left( \mathbf{E}\left| Y_t^1\right| ^2\right)
\leq C\varepsilon .
\]

The proof is complete.\quad $\Box$

\begin{lemma}
Assuming (H1)-(H4) hold, then we have
\begin{equation}
\mathbf{E}\left( \sup\limits_{0\leq t\leq T}\left| y_t^1\right| ^2\right)
\leq C\varepsilon ,  \tag{3.10}
\end{equation}
\begin{equation}
\mathbf{E}\left( \sup\limits_{0\leq t\leq T}\left| Y_t^1\right| ^2\right)
\leq C\varepsilon .  \tag{3.11}
\end{equation}
\end{lemma}

\noindent \textbf{Proof.}\quad \noindent Squaring both sides of
\begin{eqnarray*}
y_t^1 &=&\int_0^t\left( f_yy_s^1+f_YY_s^1+f_zz_s^1+f_ZZ_s^1+f\left(
u_s^\varepsilon \right) -f\left( u_s\right) \right) \text{d}s \\
&&+\int_0^t\left( g_yy_s^1+g_YY_s^1+g_zz_s^1+g_ZZ_s^1\right) \text{d}%
W_s-\int_0^tz_s^1\text{\^d}B_s,
\end{eqnarray*}
we have
\begin{eqnarray*}
\left| y_t^1\right| ^2 &\leq &3\left( \int_0^t\left(
f_yy_s^1+f_YY_s^1+f_zz_s^1+f_ZZ_s^1+f\left( u_s^\varepsilon \right) -f\left(
u_s\right) \right) \text{d}s\right) ^2 \\
&&+3\left( \int_0^t\left( g_yy_s^1+g_YY_s^1+g_zz_s^1+g_ZZ_s^1\right) \text{d}%
W_s\right) ^2+3\left( \int_0^tz_s^1\text{\^d}B_s\right) ^2 \\
&\leq &3t\int_0^t\left( f_yy_s^1+f_YY_s^1+f_zz_s^1+f_ZZ_s^1+f\left(
u_s^\varepsilon \right) -f\left( u_s\right) \right) ^2\text{d}s \\
&&+3\left( \int_0^t\left( g_yy_s^1+g_YY_s^1+g_zz_s^1+g_ZZ_s^1\right) \text{d}%
W_s\right) ^2+3\left( \int_0^Tz_s^1\text{\^d}B_s-\int_t^Tz_s^1\text{\^d}%
B_s\right) ^2 \\
&\leq &C\int_0^t\left[ \left| y_s^1\right| ^2+\left| Y_s^1\right| ^2+\left|
z_s^1\right| ^2+\left| Z_s^1\right| ^2+\left| f\left( u_s^\varepsilon
\right) -f\left( u_s\right) \right| ^2\right] \text{d}s \\
&&+3\left( \int_0^t\left( g_yy_s^1+g_YY_s^1+g_zz_s^1+g_ZZ_s^1\right) \text{d}%
W_s\right) ^2 \\
&&+6\left( \int_0^Tz_s^1\text{\^d}B_s\right) ^2+6\left( \int_t^Tz_s^1\text{%
\^d}B_s\right) ^2,
\end{eqnarray*}
then
\begin{eqnarray*}
\sup\limits_{0\leq t\leq T}\left| y_t^1\right| ^2 &\leq &C\int_0^T\left[
\left| y_s^1\right| ^2+\left| Y_s^1\right| ^2+\left| z_s^1\right| ^2+\left|
Z_s^1\right| ^2+\left| f\left( u_s^\varepsilon \right) -f\left( u_s\right)
\right| ^2\right] \text{d}s+6\left( \int_0^Tz_s^1\text{\^d}B_s\right) ^2 \\
&&\ +3\sup\limits_{0\leq t\leq T}\left( \int_0^t\left(
g_yy_s^1+g_YY_s^1+g_zz_s^1+g_ZZ_s^1\right) \text{d}W_s\right)
^2+3\sup\limits_{0\leq t\leq T}\left( \int_t^Tz_s^1\text{\^d}B_s\right) ^2 \\
\ &\leq &C\int_0^T\left[ \left| y_s^1\right| ^2+\left| Y_s^1\right|
^2+\left| z_s^1\right| ^2+\left| Z_s^1\right| ^2+\left| f\left(
u_s^\varepsilon \right) -f\left( u_s\right) \right| ^2\right] \text{d}%
s+6\left( \int_0^Tz_s^1\text{\^d}B_s\right) ^2 \\
&&\ +3\left( \sup\limits_{0\leq t\leq T}\left| \int_0^t\left(
g_yy_s^1+g_YY_s^1+g_zz_s^1+g_ZZ_s^1\right) \text{d}W_s\right| \right)
^2+3\left( \sup\limits_{0\leq t\leq T}\left| \int_t^Tz_s^1\text{\^d}%
B_s\right| \right) ^2,
\end{eqnarray*}
where $C>0$ is some constant. Hereafter, $C$ will be some generic constant,
which can be different from line to line. Taking expectation, by B-D-G
inequality and H\"older inequality, it follows that
\begin{eqnarray*}
\mathbf{E}\left( \sup\limits_{0\leq t\leq T}\left| y_t^1\right| ^2\right)
&\leq &C\mathbf{E}\int_0^T\left[ \left| y_s^1\right| ^2+\left| Y_s^1\right|
^2+\left| z_s^1\right| ^2+\left| Z_s^1\right| ^2+\left| f\left(
u_s^\varepsilon \right) -f\left( u_s\right) \right| ^2\right] \text{d}s+6%
\mathbf{E}\int_0^T\left| z_s^1\right| ^2\text{d}s \\
&&\ +C\mathbf{E}\int_0^T\left| g_yy_s^1+g_YY_s^1+g_zz_s^1+g_ZZ_s^1\right| ^2%
\text{d}s+C\mathbf{E}\int_0^T\left| z_s^1\right| ^2\text{d}s \\
\ &\leq &C\mathbf{E}\int_0^T\left[ \left| y_s^1\right| ^2+\left|
Y_s^1\right| ^2+\left| z_s^1\right| ^2+\left| Z_s^1\right| ^2\right] \text{d}%
s+C\mathbf{E}\int_0^T\left| f\left( u_s^\varepsilon \right) -f\left(
u_s\right) \right| ^2\text{d}s.
\end{eqnarray*}
From Lemma 4, (3.10) holds. By the similar argument, we can prove (3.11).
Squaring both sides of
\begin{eqnarray*}
Y_t^1 &=&h_y\left( y_T\right) y_T^1+\int_t^T\left(
F_yy_s^1+F_YY_s^1+F_zz_s^1+F_ZZ_s^1+F\left( u_s^\varepsilon \right) -F\left(
u_s\right) \right) \text{d}s \\
&&\ +\int_t^T\left( G_yy_s^1+G_YY_s^1+G_zz_s^1+G_ZZ_s^1\right) \text{\^d}%
B_s-\int_t^TZ_s^1\text{d}W_s,
\end{eqnarray*}
it follows that
\begin{eqnarray*}
\left| Y_t^1\right| ^2 &\leq &5\left| h_y\left( y_T\right) y_T^1\right|
^2+5\left( \int_t^T\left( F_yy_s^1+F_YY_s^1+F_zz_s^1+F_ZZ_s^1+F\left(
u_s^\varepsilon \right) -F\left( u_s\right) \right) \text{d}s\right) ^2 \\
&&\ +5\left( \int_t^T\left( G_yy_s^1+G_YY_s^1+G_zz_s^1+G_ZZ_s^1\right) \text{%
\^d}B_s\right) ^2+5\left( \int_0^TZ_s^1\text{d}W_s\right) ^2+5\left(
\int_0^tZ_s^1\text{d}W_s\right) ^2.
\end{eqnarray*}
Thus
\begin{eqnarray*}
\sup\limits_{0\leq t\leq T}\left| Y_t^1\right| ^2 &\leq &5\left| h_y\left(
y_T\right) y_T^1\right| ^2+5\left( \int_0^TZ_s^1\text{d}W_s\right)
^2+5\sup\limits_{0\leq t\leq T}\left( \int_0^tZ_s^1\text{d}W_s\right) ^2 \\
&&\ +5\left( T-t\right) \int_t^T\left|
F_yy_s^1+F_YY_s^1+F_zz_s^1+F_ZZ_s^1+F\left( u_s^\varepsilon \right) -F\left(
u_s\right) \right| ^2\text{d}s \\
&&\ +5\sup\limits_{0\leq t\leq T}\left( \int_t^T\left(
G_yy_s^1+G_YY_s^1+G_zz_s^1+G_ZZ_s^1\right) \text{\^d}B_s\right) ^2.
\end{eqnarray*}
Taking expectation and by B-D-G inequality, it follows that
\begin{eqnarray*}
\mathbf{E}\left( \sup\limits_{0\leq t\leq T}\left| Y_t^1\right| ^2\right)
&\leq &5\mathbf{E}\left| h_y\left( y_T\right) y_T^1\right| ^2+5\mathbf{E}%
\int_0^T\left| Z_s^1\right| ^2\text{d}s+C\mathbf{E}\int_0^T\left|
Z_s^1\right| ^2\text{d}s \\
&&\ +C\mathbf{E}\int_0^T\left( \left| y_s^1\right| ^2+\left| Y_s^1\right|
^2+\left| z_s^1\right| ^2+\left| Z_s^1\right| ^2+\left| F\left(
u_s^\varepsilon \right) -F\left( u_s\right) \right| ^2\right) \text{d}s \\
&&\ +C\mathbf{E}\int_0^T\left| G_yy_s^1+G_YY_s^1+G_zz_s^1+G_ZZ_s^1\right| ^2%
\text{d}s.
\end{eqnarray*}
Noting (3.10), from Lemma 4 and Lemma 5, it is easy to see that (3.11)
holds. The proof is complete.\quad $\Box$  Next, we will give some elaborate
estimates for $\left( y_t^1,Y_t^1,z_t^1,Z_t^1\right) $ by virtue of the
techniques of FBDSDEs.

\begin{lemma}
Assuming (H1)-(H4) hold, then we have
\begin{equation}
\mathbf{E}\int_0^T\left| y_t^1\right| ^2\text{d}t\leq C\varepsilon ^{\frac
32},  \tag{3.12}
\end{equation}
\begin{equation}
\mathbf{E}\int_0^T\left| Y_t^1\right| ^2\text{d}t\leq C\varepsilon ^{\frac
32},  \tag{3.13}
\end{equation}
\begin{equation}
\mathbf{E}\int_0^T\left| z_t^1\right| ^2\text{d}t\leq C\varepsilon ^{\frac
32},  \tag{3.14}
\end{equation}
\begin{equation}
\mathbf{E}\int_0^T\left| Z_t^1\right| ^2\text{d}t\leq C\varepsilon ^{\frac
32}.  \tag{3.15}
\end{equation}
\end{lemma}

\noindent \textbf{Proof.}\quad \noindent By (3.7), we have
\begin{eqnarray*}
&&\mathbf{E}\left[ \left| y_T^1\right| ^2h_y\left( y_T\right) \right] +\mu
\mathbf{E}\int_0^T\left( \left| y_t^1\right| ^2+\left| Y_t^1\right|
^2+\left| z_t^1\right| ^2+\left| Z_t^1\right| ^2\right) \text{d}t \\
&\leq &\mathbf{E}\int_0^T\left( f\left( u_t^\varepsilon \right) -f\left(
u_t\right) \right) Y_t^1\text{d}t-\mathbf{E}\int_0^T\left( F\left(
u_t^\varepsilon \right) -F\left( u_t\right) \right) y_t^1\text{d}t \\
&\leq &\mathbf{E}\left[ \sup_{0\leq t\leq T}\left| Y_t^1\right|
\int_0^T\left| f\left( u_t^\varepsilon \right) -f\left( u_t\right) \right|
\text{d}t\right] \\
&&+\mathbf{E}\left[ \sup_{0\leq t\leq T}\left| y_t^1\right| \text{d}%
t\int_0^T\left| F\left( u_t^\varepsilon \right) -F\left( u_t\right) \right|
\text{d}t\right] \\
&\leq &\left[ \mathbf{E}\left( \sup_{0\leq t\leq T}\left| y_t^1\right|
^2\right) \right] ^{\frac 12}\left[ \mathbf{E}\left( \int_0^T\left| F\left(
u_t^\varepsilon \right) -F\left( u_t\right) \right| \text{d}t\right)
^2\right] ^{\frac 12} \\
&&+\left[ \mathbf{E}\left( \sup_{0\leq t\leq T}\left| Y_t^1\right| ^2\right)
\right] ^{\frac 12}\left[ \mathbf{E}\left( \int_0^T\left| f\left(
u_t^\varepsilon \right) -f\left( u_t\right) \right| \text{d}t\right)
^2\right] ^{\frac 12} \\
&\leq &C\varepsilon ^{\frac 32},
\end{eqnarray*}
where $C$ is a sufficiently large positive constant. From (H3), the desired
results are obtained.\quad $\Box$

In order to obtain variational inequality, we need the following lemma.

\begin{lemma}
Assuming (H1)-(H4) hold, then we have
\begin{equation}
\mathbf{E}\int_0^T\left| y_t^\varepsilon -y_t-y_t^1\right| ^2\text{d}t\leq
C\varepsilon ^{\frac 32},  \tag{3.16}
\end{equation}
\begin{equation}
\mathbf{E}\int_0^T\left| Y_t^\varepsilon -Y_t-Y_t^1\right| ^2\text{d}t\leq
C\varepsilon ^{\frac 32},  \tag{3.17}
\end{equation}
\begin{equation}
\mathbf{E}\int_0^T\left| z_t^\varepsilon -z_t-z_t^1\right| ^2\text{d}t\leq
C\varepsilon ^{\frac 32},  \tag{3.18}
\end{equation}
\begin{equation}
\mathbf{E}\int_0^T\left| Z_t^\varepsilon -Z_t-Z_t^1\right| ^2\text{d}t\leq
C\varepsilon ^{\frac 32},  \tag{3.19}
\end{equation}
\begin{equation}
\sup\limits_{0\leq t\leq T}\left[ \mathbf{E}\left| y_t^\varepsilon
-y_t-y_t^1\right| ^2\right] \leq C\varepsilon ^{\frac 32},  \tag{3.20}
\end{equation}
\begin{equation}
\sup\limits_{0\leq t\leq T}\left[ \mathbf{E}\left| Y_t^\varepsilon
-Y_t-Y_t^1\right| ^2\right] \leq C\varepsilon ^{\frac 32}.  \tag{3.21}
\end{equation}
\end{lemma}

\noindent \textbf{Proof.}\quad \noindent For notational convenience, we
denote
\begin{eqnarray*}
\tilde y_t &=&y_t^\varepsilon -y_t-y_t^1, \\
\tilde Y_t &=&Y_t^\varepsilon -Y_t-Y_t^1, \\
\tilde z_t &=&z_t^\varepsilon -z_t-z_t^1, \\
\tilde Z_t &=&Z_t^\varepsilon -Z_t-Z_t^1.
\end{eqnarray*}
We have the following FBDSDEs
\begin{eqnarray*}
\tilde y_t &=&\int_0^t\left[ \tilde f_y\tilde y_s+\tilde f_Y\tilde
Y_s+\tilde f_z\tilde z_s+\tilde f_Z\tilde Z_s\right] \text{d}%
s+\int_0^tV_s^\varepsilon ds+\int_0^tH_s\text{d}s \\
&&+\int_0^t\left[ \tilde g_y\tilde y_s+\tilde g_Y\tilde Y_s+\tilde g_z\tilde
z_s+\tilde g_Z\tilde Z_s\right] \text{d}W_s-\int_0^t\tilde z_s\text{\^d}B_s,
\\
\tilde Y_t &=&h\left( y_T^\varepsilon \right) -h\left( y_T+y_T^1\right)
+\int_t^T\left[ \tilde F_y\tilde y_s+\tilde F_Y\tilde Y_s+\tilde F_z\tilde
z_s+\tilde F_Z\tilde Z_s\right] \text{d}s \\
&&+\int_t^T\left[ \tilde G_y\tilde y_s+\tilde G_Y\tilde Y_s+\tilde G_z\tilde
z_s+\tilde G_Z\tilde Z_s\right] \text{\^d}B_s+\int_t^T\tilde V_s^\varepsilon
\text{d}s+\int_t^T\tilde H_s\text{d}s \\
&&+\int_0^1\left( h_y\left( y_T+y_T^1\lambda \right) -h_y\left( y_T\right)
\right) y_T^1\text{d}\lambda -\int_t^T\tilde Z_s\text{d}W_s,
\end{eqnarray*}
where
\begin{eqnarray*}
\tilde f_y &=&\int_0^1f_y\left( y_s+y_s^1+\lambda \tilde
y_s,Y_s+Y_s^1+\lambda \tilde Y_s,z_s+z_s^1+\lambda \tilde
z_s,Z_s+Z_s^1+\lambda \tilde Z_s,u_s^\varepsilon \right) \text{d}\lambda , \\
\tilde f_Y &=&\int_0^1f_Y\left( y_s+y_s^1+\lambda \tilde
y_s,Y_s+Y_s^1+\lambda \tilde Y_s,z_s+z_s^1+\lambda \tilde
z_s,Z_s+Z_s^1+\lambda \tilde Z_s,u_s^\varepsilon \right) \text{d}\lambda , \\
\tilde f_z &=&\int_0^1f_z\left( y_s+y_s^1+\lambda \tilde
y_s,Y_s+Y_s^1+\lambda \tilde Y_s,z_s+z_s^1+\lambda \tilde
z_s,Z_s+Z_s^1+\lambda \tilde Z_s,u_s^\varepsilon \right) \text{d}\lambda , \\
\tilde f_Z &=&\int_0^1f_Z\left( y_s+y_s^1+\lambda \tilde
y_s,Y_s+Y_s^1+\lambda \tilde Y_s,z_s+z_s^1+\lambda \tilde
z_s,Z_s+Z_s^1+\lambda \tilde Z_s,u_s^\varepsilon \right) \text{d}\lambda ,
\end{eqnarray*}
\begin{eqnarray*}
\tilde F_y &=&\int_0^1F_y\left( y_s+y_s^1+\lambda \tilde
y_s,Y_s+Y_s^1+\lambda \tilde Y_s,z_s+z_s^1+\lambda \tilde
z_s,Z_s+Z_s^1+\lambda \tilde Z_s,u_s^\varepsilon \right) \text{d}\lambda , \\
\tilde F_Y &=&\int_0^1F_Y\left( y_s+y_s^1+\lambda \tilde
y_s,Y_s+Y_s^1+\lambda \tilde Y_s,z_s+z_s^1+\lambda \tilde
z_s,Z_s+Z_s^1+\lambda \tilde Z_s,u_s^\varepsilon \right) \text{d}\lambda , \\
\tilde F_z &=&\int_0^1F_z\left( y_s+y_s^1+\lambda \tilde
y_s,Y_s+Y_s^1+\lambda \tilde Y_s,z_s+z_s^1+\lambda \tilde
z_s,Z_s+Z_s^1+\lambda \tilde Z_s,u_s^\varepsilon \right) \text{d}\lambda , \\
\tilde F_Z &=&\int_0^1F_Z\left( y_s+y_s^1+\lambda \tilde
y_s,Y_s+Y_s^1+\lambda \tilde Y_s,z_s+z_s^1+\lambda \tilde
z_s,Z_s+Z_s^1+\lambda \tilde Z_s,u_s^\varepsilon \right) \text{d}\lambda ,
\end{eqnarray*}
\begin{eqnarray*}
\tilde g_y &=&\int_0^1g_y\left( y_s+y_s^1+\lambda \tilde
y_s,Y_s+Y_s^1+\lambda \tilde Y_s,z_s+z_s^1+\lambda \tilde
z_s,Z_s+Z_s^1+\lambda \tilde Z_s\right) \text{d}\lambda , \\
\tilde g_Y &=&\int_0^1g_Y\left( y_s+y_s^1+\lambda \tilde
y_s,Y_s+Y_s^1+\lambda \tilde Y_s,z_s+z_s^1+\lambda \tilde
z_s,Z_s+Z_s^1+\lambda \tilde Z_s\right) \text{d}\lambda , \\
\tilde g_z &=&\int_0^1g_z\left( y_s+y_s^1+\lambda \tilde
y_s,Y_s+Y_s^1+\lambda \tilde Y_s,z_s+z_s^1+\lambda \tilde
z_s,Z_s+Z_s^1+\lambda \tilde Z_s\right) \text{d}\lambda , \\
\tilde g_Z &=&\int_0^1g_Z\left( y_s+y_s^1+\lambda \tilde
y_s,Y_s+Y_s^1+\lambda \tilde Y_s,z_s+z_s^1+\lambda \tilde
z_s,Z_s+Z_s^1+\lambda \tilde Z_s\right) \text{d}\lambda ,
\end{eqnarray*}
\begin{eqnarray*}
\tilde G_y &=&\int_0^1G_y\left( y_s+y_s^1+\lambda \tilde
y_s,Y_s+Y_s^1+\lambda \tilde Y_s,z_s+z_s^1+\lambda \tilde
z_s,Z_s+Z_s^1+\lambda \tilde Z_s\right) \text{d}\lambda , \\
\tilde G_Y &=&\int_0^1G_Y\left( y_s+y_s^1+\lambda \tilde
y_s,Y_s+Y_s^1+\lambda \tilde Y_s,z_s+z_s^1+\lambda \tilde
z_s,Z_s+Z_s^1+\lambda \tilde Z_s\right) \text{d}\lambda , \\
\tilde G_z &=&\int_0^1G_z\left( y_s+y_s^1+\lambda \tilde
y_s,Y_s+Y_s^1+\lambda \tilde Y_s,z_s+z_s^1+\lambda \tilde
z_s,Z_s+Z_s^1+\lambda \tilde Z_s\right) \text{d}\lambda , \\
\tilde G_Z &=&\int_0^1G_Z\left( y_s+y_s^1+\lambda \tilde
y_s,Y_s+Y_s^1+\lambda \tilde Y_s,z_s+z_s^1+\lambda \tilde
z_s,Z_s+Z_s^1+\lambda \tilde Z_s\right) \text{d}\lambda ,
\end{eqnarray*}
\begin{eqnarray*}
V_s^\varepsilon &=&\int_0^1\left[ f_y\left( y_s+\lambda y_s^1,Y_s+\lambda
Y_s^1,z_s+\lambda z_s^1,Z_s+\lambda Z_s^1,u_s^\varepsilon \right)
-f_y\right] y_s^1\text{d}\lambda \\
&&+\int_0^1\left[ f_Y\left( y_s+\lambda y_s^1,Y_s+\lambda Y_s^1,z_s+\lambda
z_s^1,Z_s+\lambda Z_s^1,u_s^\varepsilon \right) -f_Y\right] Y_s^1\text{d}%
\lambda \\
&&+\int_0^1\left[ f_z\left( y_s+\lambda y_s^1,Y_s+\lambda Y_s^1,z_s+\lambda
z_s^1,Z_s+\lambda Z_s^1,u_s^\varepsilon \right) -f_z\right] z_s^1\text{d}%
\lambda \\
&&+\int_0^1\left[ f_Z\left( y_s+\lambda y_s^1,Y_s+\lambda Y_s^1,z_s+\lambda
z_s^1,Z_s+\lambda Z_s^1,u_s^\varepsilon \right) -f_Z\right] Z_s^1\text{d}%
\lambda ,
\end{eqnarray*}
\begin{eqnarray*}
H_s &=&\int_0^1\left[ g_y\left( y_s+y_s^1\lambda ,Y_s+\lambda
Y_s^1,z_s+\lambda z_s^1,Z_s+\lambda Z_s^1\right) -g_y\right] y_s^1\text{d}%
\lambda \\
&&+\int_0^1\left[ g_Y\left( y_s+y_s^1\lambda ,Y_s+\lambda Y_s^1,z_s+\lambda
z_s^1,Z_s+\lambda Z_s^1\right) -g_Y\right] Y_s^1\text{d}\lambda \\
&&+\int_0^1\left[ g_z\left( y_s+y_s^1\lambda ,Y_s+\lambda Y_s^1,z_s+\lambda
z_s^1,Z_s+\lambda Z_s^1\right) -g_z\right] z_s^1\text{d}\lambda \\
&&+\int_0^1\left[ g_Z\left( y_s+y_s^1\lambda ,Y_s+\lambda Y_s^1,z_s+\lambda
z_s^1,Z_s+\lambda Z_s^1\right) -g_Z\right] Z_s^1\text{d}\lambda ,
\end{eqnarray*}
\begin{eqnarray*}
\tilde V_s^\varepsilon &=&\int_0^1\left[ F_y\left( y_s+y_s^1\lambda
,Y_s+\lambda Y_s^1,z_s+\lambda z_s^1,Z_s+\lambda Z_s^1,u_s^\varepsilon
\right) -F_y\right] y_s^1\text{d}\lambda \\
&&+\int_0^1\left[ F_Y\left( y_s+y_s^1\lambda ,Y_s+\lambda Y_s^1,z_s+\lambda
z_s^1,Z_s+\lambda Z_s^1,u_s^\varepsilon \right) -F_Y\right] Y_s^1\text{d}%
\lambda \\
&&+\int_0^1\left[ F_z\left( y_s+y_s^1\lambda ,Y_s+\lambda Y_s^1,z_s+\lambda
z_s^1,Z_s+\lambda Z_s^1,u_s^\varepsilon \right) -F_z\right] z_s^1\text{d}%
\lambda \\
&&+\int_0^1\left[ F_Z\left( y_s+y_s^1\lambda ,Y_s+\lambda Y_s^1,z_s+\lambda
z_s^1,Z_s+\lambda Z_s^1,u_s^\varepsilon \right) -F_Z\right] Z_s^1\text{d}%
\lambda ,
\end{eqnarray*}
\begin{eqnarray*}
\tilde H_s &=&\int_0^1\left[ G_y\left( y_s+y_s^1\lambda ,Y_s+\lambda
Y_s^1,z_s+\lambda z_s^1,Z_s+\lambda Z_s^1\right) -G_y\right] y_s^1\text{d}%
\lambda \\
&&+\int_0^1\left[ G_Y\left( y_s+y_s^1\lambda ,Y_s+\lambda Y_s^1,z_s+\lambda
z_s^1,Z_s+\lambda Z_s^1\right) -G_Y\right] Y_s^1\text{d}\lambda \\
&&+\int_0^1\left[ G_z\left( y_s+y_s^1\lambda ,Y_s+\lambda Y_s^1,z_s+\lambda
z_s^1,Z_s+\lambda Z_s^1\right) -G_z\right] z_s^1\text{d}\lambda \\
&&+\int_0^1\left[ G_Z\left( y_s+y_s^1\lambda ,Y_s+\lambda Y_s^1,z_s+\lambda
z_s^1,Z_s+\lambda Z_s^1\right) -G_Z\right] Z_s^1\text{d}\lambda .
\end{eqnarray*}
It is easy to check that
\begin{eqnarray*}
\mathbf{E}\int_0^T\left| \tilde V_s^\varepsilon \right| ^2\text{d}s &\leq
&C\varepsilon ^{\frac 32}, \\
\mathbf{E}\int_0^T\left| \tilde H_s\right| ^2\text{d}s &\leq &C\varepsilon
^{\frac 32}, \\
\mathbf{E}\int_0^T\left| V_s^\varepsilon \right| ^2\text{d}s &\leq
&C\varepsilon ^{\frac 32}, \\
\mathbf{E}\int_0^T\left| H_s\right| ^2\text{d}s &\leq &C\varepsilon ^{\frac
32}.
\end{eqnarray*}
By Lemma 7, applying It\^o's formula to $\left\langle \tilde y_t,\tilde
Y_t\right\rangle $ on $\left[ 0,T\right] ,$ we get
\begin{eqnarray*}
&&\mathbf{E}\left\langle h\left( y_T^\varepsilon \right) -h\left( y_T\right)
-h_y\left( y_T\right) y_T^1,y_T^\varepsilon -y_T-y_T^1\right\rangle \\
&&+\mu \mathbf{E}\int_0^T\left[ \left| \tilde y_s\right| ^2+\left| \tilde
Y_s\right| ^2+\left| \tilde z_s\right| ^2+\left| \tilde Z_s\right| ^2\right]
\text{d}s \\
&\leq &-\mathbf{E}\int_0^T\left\langle \tilde y_s,\tilde V_s^\varepsilon
+\tilde H_s\right\rangle \text{d}s+\mathbf{E}\int_0^T\left\langle \tilde
Y_s,V_s^\varepsilon +H_s\right\rangle \text{d}s \\
&\leq &\mathbf{E}\frac \mu 2\int_0^T\left| \tilde y_s\right| ^2\text{d}s+%
\mathbf{E}\frac 1\mu \int_0^T\left| \tilde V_s^\varepsilon \right| ^2\text{d}%
s+\frac 1\mu \mathbf{E}\int_0^T\left| \tilde H_s\right| ^2\text{d}s \\
&&+\mathbf{E}\frac \mu 2\int_0^T\left| \tilde Y_s\right| ^2\text{d}s+\mathbf{%
E}\frac 1\mu \int_0^T\left| V_s^\varepsilon \right| ^2\text{d}s+\frac 1\mu
\mathbf{E}\int_0^T\left| H_s\right| ^2\text{d}s.
\end{eqnarray*}
Noting that by means of the same arguments in Lemma 5, from Lemma 7, we
easily have
\begin{eqnarray*}
\sup\limits_{0\leq t\leq T}\left( \mathbf{E}\left| y_t^1\right| ^2\right)
&\leq &C\varepsilon ^{\frac 32}, \\
\sup\limits_{0\leq t\leq T}\left( \mathbf{E}\left| Y_t^1\right| ^2\right)
&\leq &C\varepsilon ^{\frac 32}.
\end{eqnarray*}
Thus it is obvious that
\[
\mathbf{E}h\left( y_T+y_T^1\right) =\mathbf{E}h\left( y_T\right) +\mathbf{E}%
h_y\left( y_T\right) y_T^1+C\varepsilon ^{\frac 32},
\]
so by (H3) it follows that
\begin{eqnarray*}
&&\mathbf{E}\left\langle h\left( y_T^\varepsilon \right) -h\left( y_T\right)
-h_y\left( y_T\right) y_T^1,y_T^\varepsilon -y_T-y_T^1\right\rangle \\
&=&\mathbf{E}\left\langle h\left( y_T^\varepsilon \right) -h\left(
y_T+y_T^1\right) +C\varepsilon ^{\frac 32},y_T^\varepsilon
-y_T-y_T^1\right\rangle \\
&\geq &\mathbf{E}\left( y_T^\varepsilon -y_T-y_T^1\right) \cdot C\varepsilon
^{\frac 32},
\end{eqnarray*}
and
\begin{eqnarray*}
&&\mathbf{E}\mu \int_0^T\left[ \frac 12\left| \tilde y_s\right| ^2+\frac
12\left| \tilde Y_s\right| ^2+\left| \tilde z_s\right| ^2+\left| \tilde
Z_s\right| ^2\right] \text{d}s \\
&\leq &\mathbf{E}\frac 1\mu \int_0^T\left| \tilde V_s^\varepsilon \right| ^2%
\text{d}s+\frac 1\mu \mathbf{E}\int_0^T\left| \tilde H_s\right| ^2\text{d}s
\\
&&+\mathbf{E}\frac 1\mu \int_0^T\left| V_s^\varepsilon \right| ^2\text{d}%
s+\frac 1\mu \mathbf{E}\int_0^T\left| H_s\right| ^2\text{d}s-\mathbf{E}%
\left( \tilde y_T\right) \cdot C\varepsilon ^{\frac 32}.
\end{eqnarray*}
It is not difficult to see that $\mathbf{E}\left( \tilde y_T\right) $ is
bounded. Consequently, from that, (3.16)-(3.19) hold. Further using the
similar arguments in Lemma 5, we can obtain (3.20) and (3.21). The proof is
complete.\quad $\Box$

\begin{lemma}
\textbf{(Variational inequality) }Under the assumptions (H1)-(H4), it holds
that
\begin{equation}
\mathbf{E}\int_0^T\left[ l_yy_t^1+l_YY_t^1+l_zz_t^1+l_ZZ_t^1+l\left(
u_t^\varepsilon \right) -l\left( u_t\right) \right] \text{d}t+\mathbf{E}%
\left[ \Phi _y\left( y_T\right) y_T^1\right] +\mathbf{E}\left[ \gamma
_Y\left( Y_0\right) Y_0^1\right] \geq o\left( \varepsilon \right) .
\tag{3.22}
\end{equation}
\end{lemma}

\noindent \textbf{Proof.}\quad \noindent According to the definition of $%
u_t^\varepsilon ,$ we have
\[
J\left( u_{\left( \cdot \right) }^\varepsilon \right) \geq J\left( u_{\left(
\cdot \right) }\right) ,
\]
moreover
\begin{eqnarray*}
&&\mathbf{E}\int_0^T\left[ l\left( t,y_t^\varepsilon ,Y_t^\varepsilon
,z_t^\varepsilon ,Z_t^\varepsilon ,u_t^\varepsilon \right) -l\left(
t,y_t,Y_t,z_t,Z_t,u_t\right) \right] \text{d}t \\
&&+\mathbf{E}\left[ \Phi \left( y_T^\varepsilon \right) -\Phi \left(
y_T\right) \right] +\mathbf{E}\left[ \gamma \left( Y_0^\varepsilon \right)
-\gamma \left( Y_0\right) \right] \\
&\geq &0,
\end{eqnarray*}
or
\begin{eqnarray*}
&&\mathbf{E}\int_0^T\left[ l\left( t,y_t^\varepsilon ,Y_t^\varepsilon
,z_t^\varepsilon ,Z_t^\varepsilon ,u_t^\varepsilon \right) -l\left(
t,y_t+y_t^1,Y_t+Y_t^1,z_t+z_t^1,Z_t+Z_t^1,u_t^\varepsilon \right) \right]
\text{d}t \\
&&+\mathbf{E}\int_0^T\left[ l\left(
t,y_t+y_t^1,Y_t+Y_t^1,z_t+z_t^1,Z_t+Z_t^1,u_t^\varepsilon \right) -l\left(
t,y_t,Y_t,z_t,Z_t,u_t\right) \right] \text{d}t \\
&&+\mathbf{E}\left[ \Phi \left( y_T^\varepsilon \right) -\Phi \left(
y_T+y_T^1\right) \right] +\mathbf{E}\left[ \Phi \left( y_T+y_T^1\right)
-\Phi \left( y_T\right) \right] \\
&&+\mathbf{E}\left[ \gamma \left( Y_0^\varepsilon \right) -\gamma \left(
Y_0+Y_T^1\right) \right] +\mathbf{E}\left[ \gamma \left( Y_0+Y_T^1\right)
-\gamma \left( Y_0\right) \right] \\
&\geq &0.
\end{eqnarray*}
By Lemma 8, it follows that
\begin{eqnarray*}
&&\mathbf{E}\int_0^T\left[ l\left( t,y_t^\varepsilon ,Y_t^\varepsilon
,z_t^\varepsilon ,Z_t^\varepsilon ,u_t^\varepsilon \right) -l\left(
t,y_t+y_t^1,Y_t+Y_t^1,z_t+z_t^1,Z_t+Z_t^1,u_t^\varepsilon \right) \right]
\text{d}t \\
&&+\mathbf{E}\left[ \Phi \left( y_T^\varepsilon \right) -\Phi \left(
y_T+y_T^1\right) \right] +\mathbf{E}\left[ \gamma \left( Y_0^\varepsilon
\right) -\gamma \left( Y_0+Y_T^1\right) \right] \\
&\leq &C\varepsilon ^{\frac 32},
\end{eqnarray*}
while
\begin{eqnarray*}
0 &\leq &\mathbf{E}\int_0^T\left[ l\left(
t,y_t+y_t^1,Y_t+Y_t^1,z_t+z_t^1,Z_t+Z_t^1,u_t^\varepsilon \right) -l\left(
t,y_t,Y_t,z_t,Z_t,u_t\right) \right] \text{d}t \\
&&\ +\mathbf{E}\left[ \Phi \left( y_T+y_T^1\right) -\Phi \left( y_T\right)
\right] +\mathbf{E}\left[ \gamma \left( Y_0+Y_T^1\right) -\gamma \left(
Y_0\right) \right] +C\varepsilon ^{\frac 32} \\
\ &=&\mathbf{E}\int_0^T\left[ l\left(
t,y_t+y_t^1,Y_t+Y_t^1,z_t+z_t^1,Z_t+Z_t^1,u_t\right) -l\left(
t,y_t,Y_t,z_t,Z_t,u_t\right) \right] \text{d}t \\
&&\ +\mathbf{E}\int_0^Tl\left(
t,y_t+y_t^1,Y_t+Y_t^1,z_t+z_t^1,Z_t+Z_t^1,u_t^\varepsilon \right) \text{d}t
\\
&&\ -\mathbf{E}\int_0^Tl\left(
t,y_t+y_t^1,Y_t+Y_t^1,z_t+z_t^1,Z_t+Z_t^1,u_t\right) \text{d}t \\
&&\ +\mathbf{E}\left[ \Phi \left( y_T+y_T^1\right) -\Phi \left( y_T\right)
\right] +\mathbf{E}\left[ \gamma \left( Y_0+Y_T^1\right) -\gamma \left(
Y_0\right) \right] +C\varepsilon ^{\frac 32} \\
\ &=&\mathbf{E}\int_0^T\left[ l_yy_t^1+l_YY_t^1+l_zz_t^1+l_ZZ_t^1\right]
\text{d}t+\mathbf{E}\int_0^T\left[ l\left( u_t^\varepsilon \right) -l\left(
u_t\right) \right] \text{d}t \\
&&\ +\mathbf{E}\int_0^T\left\{ \left[ l_y\left( u_t^\varepsilon \right)
-l_y\left( u_t\right) \right] y_t^1+\left[ l_Y\left( u_t^\varepsilon \right)
-l_Y\left( u_t\right) \right] Y_t^1\right\} \text{d}t \\
&&\ +\mathbf{E}\int_0^T\left\{ \left[ l_z\left( u_t^\varepsilon \right)
-l_z\left( u_t\right) \right] z_t^1+\left[ l_Z\left( u_t^\varepsilon \right)
-l_Z\left( u_t\right) \right] Z_t^1\right\} \text{d}t \\
&&\ +\mathbf{E}\left[ \Phi _y\left( y_T\right) y_T^1\right] +\mathbf{E}%
\left[ \gamma _Y\left( Y_0\right) Y_0^1\right] +C\varepsilon ^{\frac 32} \\
\ &\leq &\mathbf{E}\int_0^T\left[ l_yy_t^1+l_YY_t^1+l_zz_t^1+l_ZZ_t^1\right]
\text{d}t+\mathbf{E}\int_0^T\left[ l\left( u_t^\varepsilon \right) -l\left(
u_t\right) \right] \text{d}t \\
&&\ +\mathbf{E}\left[ \sup\limits_{0\leq t\leq T}\left| y_t^1\right|
\int_0^T\left| l_y\left( u_t^\varepsilon \right) -l_y\left( u_t\right)
\right| \text{d}t\right] \\
&&\ +\mathbf{E}\left[ \sup\limits_{0\leq t\leq T}\left| Y_t^1\right| \text{d}%
t\int_0^T\left| l_Y\left( u_t^\varepsilon \right) -l_Y\left( u_t\right)
\right| \text{d}t\right] \\
&&\ +\left[ \mathbf{E}\int_0^T\left| l_z\left( u_t^\varepsilon \right)
-l_z\left( u_t\right) \right| ^2\text{d}t\right] ^{\frac 12}\left[ \mathbf{E}%
\int_0^T\left| z_t^1\right| ^2\text{d}t\right] ^{\frac 12} \\
&&\ +\left[ \mathbf{E}\int_0^T\left| l_Z\left( u_t^\varepsilon \right)
-l_Z\left( u_t\right) \right| ^2\text{d}t\right] ^{\frac 12}\left[ \mathbf{E}%
\int_0^T\left| Z_t^1\right| ^2\text{d}t\right] ^{\frac 12} \\
&&\ \ \ +\mathbf{E}\left[ \Phi _y\left( y_T\right) y_T^1\right] +\mathbf{E}%
\left[ \gamma _Y\left( Y_0\right) Y_0^1\right] +C\varepsilon ^{\frac 32} \\
\ &\leq &\mathbf{E}\int_0^T\left[ l_yy_t^1+l_YY_t^1+l_zz_t^1+l_ZZ_t^1\right]
\text{d}t+\mathbf{E}\int_0^T\left[ l\left( u_t^\varepsilon \right) -l\left(
u_t\right) \right] \text{d}t \\
&&+\left[ \mathbf{E}\left( \sup\limits_{0\leq t\leq T}\left| y_t^1\right|
^2\right) \right] ^{\frac 12}\left[ \mathbf{E}\left( \int_0^T\left|
l_y\left( u_t^\varepsilon \right) -l_y\left( u_t\right) \right| \text{d}%
t\right) ^2\right] ^{\frac 12} \\
&&+\left[ \mathbf{E}\sup\limits_{0\leq t\leq T}\left( \left| Y_t^1\right|
^2\right) \right] ^{\frac 12}\left[ \mathbf{E}\left( \int_0^T\left|
l_Y\left( u_t^\varepsilon \right) -l_Y\left( u_t\right) \right| \text{d}%
t\right) ^2\right] ^{\frac 12} \\
&&\ \ \  \\
&&\ \ \ \ \ +C\varepsilon ^{\frac 12}\cdot C\varepsilon ^{\frac
34}+C\varepsilon ^{\frac 12}\cdot C\varepsilon ^{\frac 34}+\mathbf{E}\left[
\Phi _y\left( y_T\right) y_T^1\right] +\mathbf{E}\left[ \gamma _Y\left(
Y_0\right) Y_0^1\right] +C\varepsilon ^{\frac 32} \\
\ &=&\mathbf{E}\int_0^T\left[ l_yy_t^1+l_YY_t^1+l_zz_t^1+l_ZZ_t^1+l\left(
u_t^\varepsilon \right) -l\left( u_t\right) \right] dt \\
&&\ \ \ \ \ \ \ \ \ \ +\mathbf{E}\left[ \Phi _y\left( y_T\right)
y_T^1\right] +\mathbf{E}\left[ \gamma _Y\left( Y_0\right) Y_0^1\right]
+o\left( \varepsilon \right) .
\end{eqnarray*}
From that, the desired result is obtained. \quad $\Box$

\section{The maximum principle in global form}

\hspace{0.25in}We introduce the adjoint equations by virtue of dual
technique and Hamilton function for our control problem. From the
variational inequality obtained in Lemma 9, the maximum principle can be
proved by means of It\^o's formula. The adjoint equations are as follows:
\begin{equation}
\left\{
\begin{array}{l}
\text{d}p_t=(F_Yp_t-f_Yq_t+G_Yk_t-g_Yh_t-l_Y)\text{d}t \\
\quad \quad +(F_Zp_t-f_Zq_t+G_Zk_t-g_Zh_t-l_Z)\text{d}W_t-k_t\text{\^d}B_t,
\\
\text{d}q_t=(F_yp_t-f_yq_t+G_yk_t-g_yh_t-l_y)\text{d}t \\
\qquad +(F_zp_t-f_zq_t+G_zk_t-g_zh_t-l_z)\text{\^d}B_t+h_t\text{d}W_t, \\
p_0=-\gamma _Y\left( Y_0\right) ,\quad \quad q_T=-h_y\left( y_T\right)
P_T+\Phi _y\left( y_T\right) ,\quad \quad 0\leq t\leq T,
\end{array}
\right.  \tag{4.1}
\end{equation}
where $\left( p_{\left( \cdot \right) },q_{\left( \cdot \right) },k_{\left(
\cdot \right) },h_{\left( \cdot \right) }\right) \in \mathbf{R\times R\times
R}^l\mathbf{\times R}^d.$ It is easy to verify that FBDSDE (4.1) satisfies
(H1), (H2) and (H3)'. From Proposition 2, we know that (4.1) has a unique
solution $\left( p_{\left( \cdot \right) },q_{\left( \cdot \right)
},k_{\left( \cdot \right) },h_{\left( \cdot \right) }\right) \in M^2\left(
0,T;\mathbf{R\times R\times R}^l\mathbf{\times R}^d\right) .$ Now we define
the Hamilton function as follows:
\begin{eqnarray*}
H\left( t,y,Y,z,Z,v,p,q,k,h\right) &\doteq &\left\langle q,f\left(
t,y,Y,z,Z,v\right) \right\rangle -\left\langle p,F\left( t,y,Y,z,Z,v\right)
\right\rangle \\
&&\ -\left\langle k,G\left( t,y,Y,z,Z\right) \right\rangle +\left\langle
h,g\left( t,y,Y,z,Z\right) \right\rangle \\
&&\ +l\left( t,y,Y,z,Z,v\right) ,
\end{eqnarray*}
\begin{equation}
\tag{4.2}
\end{equation}

where $H:\left[ 0,T\right] \mathbf{\times R\times R\times R}^l\mathbf{\times
R}^d\mathbf{\times R\times R\times R\times R}^l\mathbf{\times R}^d\mathbf{%
\rightarrow R.}$ (4.1) can be rewritten as
\begin{equation}
\left\{
\begin{array}{l}
\text{d}p_t=-H_Y\text{d}t-H_Z\text{d}W_t-k_t\text{\^d}B_t, \\
\text{d}q_t=-H_y\text{d}t-H_z\text{\^d}B_t+h_t\text{d}W_t, \\
p_0=-\gamma _Y\left( Y_0\right) , \\
q_T=-h_y\left( y_T\right) P_T+\Phi _y\left( y_T\right) ,\quad 0\leq t\leq T.
\end{array}
\right.  \tag{4.3}
\end{equation}

From Lemma 9 and (4.2), we can obtain the main result in this paper.

\begin{theorem}
Suppose (H1)-(H4) hold. Let $\left( y_{\left( \cdot \right) },Y_{\left(
\cdot \right) },z_{\left( \cdot \right) },Z_{\left( \cdot \right)
},u_{\left( \cdot \right) }\right) $ be an optimal control and its
corresponding trajectory of (2.1), $\left( p_{\left( \cdot \right)
},q_{\left( \cdot \right) },k_{\left( \cdot \right) },h_{\left( \cdot
\right) }\right) $ be the corresponding solution of (4.1). Then the maximum
principle holds, that is
\begin{eqnarray*}
&&\ H\left( t,y_t,Y_t,z_t,Z_t,v,p_t,q_t,k_t,h_t\right) \\
\ &\geq &H\left( t,y_t,Y_t,z_t,Z_t,u_t,p_t,q_t,k_t,h_t\right) , \\
\forall v &\in &\mathcal{U},\text{ a.e, a.s..}
\end{eqnarray*}
\begin{equation}
\tag{4.4}
\end{equation}
\end{theorem}

\noindent \textbf{Proof.}\quad\noindent By applying It\^o's formula to $%
\left\langle p_t,Y_t^1\right\rangle +\left\langle q_t,y_t^1\right\rangle $,
and noting the variational equation (3.1), the adjoint equation (4.1) and
the variational inequality (3.22), we get
\begin{eqnarray*}
&&\quad \mathbf{E}\left[ \Phi _y\left( y_T\right) y_T^1\right] +\mathbf{E}%
\left[ \gamma _Y\left( Y_0\right) Y_0^1\right] \\
&&+\mathbf{E}\int_0^T\left[ l_yy_t^1+l_YY_t^1+l_zz_t^1+l_ZZ_t^1+l\left(
u_t^\varepsilon \right) -l\left( u_t\right) \right] \text{d}t \\
&=&\mathbf{E}\int_0^T\left[ H\left( t,y_t,Y_t,z_t,Z_t,u_t^\varepsilon
,p_t,q_t,k_t,h_t\right) -H\left(
t,y_t,Y_t,z_t,Z_t,u_t,p_t,q_t,k_t,h_t\right) \right] \text{d}t \\
&\geq &o\left( \varepsilon \right) .
\end{eqnarray*}
Since $\varepsilon >0$ can be arbitrarily small, from the above inequality,
(4.4) can be easily obtained. The proof is complete. \quad $\Box$

In the last part of this section, we provide a concrete example of
forward-backward doubly stochastic LQ control problems. We give the explicit
optimal control and validate our major theoretical results in Theorem 10.

\begin{example}
Let the control domain be $\mathcal{U}=\left[ -1,1\right] .$ Consider the
following linear forward-backward doubly stochastic control system which is
a simple case of (2.1). We assume that $l=d=1.$%
\begin{equation}
\left\{
\begin{array}{l}
\text{d}y_t=\left( z_t-Z_t+v_t\right) \text{d}W_t-z_t\text{\^d}B_t, \\
\text{d}Y_t=-\left( z_t+Z_t+v_t\right) \text{\^d}B_t+Z_t\text{d}W_t, \\
y_0=0,\quad Y_T=0,\text{ }\quad t\in \left[ 0,T\right] ,
\end{array}
\right.  \tag{4.5}
\end{equation}
where $T>0$ is a given constant and the cost function is
\begin{equation}
J\left( v_{\left( \cdot \right) }\right) =\frac 12\mathbf{E}\int_0^T\left(
y_t^2+Y_t^2+z_t^2+Z_t^2+v_t^2\right) \text{d}t+\frac 12\mathbf{E}y_T^2+\frac
12\mathbf{E}Y_0^2.  \tag{4.6}
\end{equation}
\end{example}

Note that (4.5) are a linear control system. According to the existence and
uniqueness for (4.5), it is straightforward to know the optimal control is $%
u_{\left( \cdot \right) }\equiv 0,$ with the corresponding optimal state
trajectory $\left( y_t,Y_t,z_t,Z_t\right) \equiv 0,$ $t\in \left[ 0,T\right]
.$ Notice that the adjoint equation associated with the optimal quadruple $%
\left( y_t,Y_t,z_t,Z_t\right) \equiv 0$ are
\begin{equation}
\left\{
\begin{array}{l}
\text{d}p_t=-Y_t\text{d}t+\left( k_t-Z_t\right) \text{d}W_t-k_t\text{\^d}B_t,
\\
\text{d}q_t=-y_t\text{d}t+\left( k_t-z_t\right) \text{\^d}B_t+h_t\text{d}W_t,
\\
p_0=0,\quad q_T=0,\quad t\in \left[ 0,T\right] .
\end{array}
\right.  \tag{4.7}
\end{equation}
Obviously, $\left( p_t,q_t,k_t,h_t\right) \equiv 0$ is the unique solution
of (4.7). Instantly, we give the Hamilton function is
\begin{eqnarray*}
H\left( t,y_t,Y_t,z_t,Z_t,v,p_t,q_t,k_t,h_t\right) &=&\frac 12\left(
y_t^2+Y_t^2+z_t^2+Z_t^2+v^2\right) \\
&&\ \ -k_t\left( z_t+Z_t+v\right) \\
&&\ \ +h_t\left( z_t-Z_t+v\right) \\
\ &=&\frac 12v^2.
\end{eqnarray*}
It is clear that, for any $v\in \mathcal{U}$, we always have
\[
H\left( t,y_t,Y_t,z_t,Z_t,v,p_t,q_t,k_t,h_t\right) \geq H\left(
t,y_t,Y_t,z_t,Z_t,u_t,p_t,q_t,k_t,h_t\right) =0,\quad \text{ a.e, a.s.}.
\]

\section{Applications to optimal control problems of stochastic partial
differential equations}

\hspace{0.25in}Let us first give some notations from [14]. For convenience,
all the variables in this section are one-dimensional. From now on $%
C^k\left( \mathbf{R};\mathbf{R}\right) ,$ $C_{l,b}^k\left( \mathbf{R};%
\mathbf{R}\right) , $ $C_p^k\left( \mathbf{R};\mathbf{R}\right) $ will
denote respectively the set of functions of class $C^k$ from $\mathbf{R}$
into $\mathbf{R}$, the set of those functions of class $C^k$ whose partial
derivatives of order less than or equal to $k$ are bounded (and hence the
function itself grows at most linearly at infinity), and the set of those
functions of class $C^k$ which, together with all their partial derivatives
of order less than or equal to $k $, grow at most like a polynomial function
of the variable $x$ at infinity. We consider the following quasilinear SPDEs
with control variable:
\begin{equation}
\left\{
\begin{array}{c}
u\left( t,x\right) =\tilde h\left( x\right) +\int_t^T\left[ \mathcal{L}%
u\left( s,x\right) +f\left( s,x,u\left( s,x\right) ,\left( \nabla u\sigma
\right) \left( s,x\right) ,v_s\right) \right] \text{d}s \\
+\int_t^Tg\left( s,x,u\left( s,x\right) ,\left( \nabla u\sigma \right)
\left( s,x\right) \right) \text{\^d}B_s,\quad 0\leq t\leq T,
\end{array}
\right.  \tag{5.1}
\end{equation}
where $u:\left[ 0,T\right] \times \mathbf{R}\rightarrow \mathbf{R}$ and $%
\nabla u\left( s,x\right) $ denotes the first order derivative of $u\left(
s,x\right) $ with respect to $x$, and
\[
\mathcal{L}u=\left(
\begin{array}{c}
Lu_1 \\
\vdots \\
Lu_k
\end{array}
\right) ,
\]
with $L\phi \left( x\right) =\frac 12\sum_{i,j=1}^d\left( \sigma \sigma
^{*}\right) _{ij}\left( x\right) \frac{\partial ^2\phi \left( x\right) }{%
\partial x_i\partial x_j}+\sum_{i=1}^db_i\left( x,v\right) \frac{\partial
\phi \left( x\right) }{\partial x_i}.$ In the present paper, we set $d=k=1,$
and
\begin{eqnarray*}
b &:&\mathbf{R\times R}\rightarrow \mathbf{R,} \\
\sigma &:&\mathbf{R}\rightarrow \mathbf{R,} \\
f &:&\left[ 0,T\right] \times \mathbf{R\times R\times R\times R}\rightarrow
\mathbf{R,} \\
g &:&\left[ 0,T\right] \times \mathbf{R\times R\times R}\rightarrow \mathbf{%
R,} \\
\tilde h &:&\mathbf{R}\rightarrow \mathbf{R}.
\end{eqnarray*}

In order to assure the existence and uniqueness of solutions for (5.1) and
(5.3) below, we give the following assumptions for sake of completeness (see
[14] for more details).

\begin{enumerate}
\item[(A1)]
\[
\left\{
\begin{array}{l}
b\in C_{l,b}^3\left( \mathbf{R\times R};\mathbf{R}\right) ,\quad \sigma \in
C_{l,b}^3\left( \mathbf{R};\mathbf{R}\right) ,\text{\quad }\tilde h\in
C_p^3\left( \mathbf{R};\mathbf{R}\right) , \\
f\left( t,\cdot ,\cdot ,\cdot ,v\right) \in C_{l,b}^3\left( \mathbf{R\times
R\times R};\mathbf{R}\right) ,\quad f\left( \cdot ,x,y,z,v\right) \in
M^2\left( 0,T;\mathbf{R}\right) , \\
g\left( t,\cdot ,\cdot ,\cdot \right) \in C_{l,b}^3\left( \mathbf{R\times
R\times R};\mathbf{R}\right) ,\quad g\left( \cdot ,x,y,z\right) \in
M^2\left( 0,T;\mathbf{R}\right) \\
\forall t\in \left[ 0,T\right] \text{, }x\in \mathbf{R}\text{, }y\in \mathbf{%
R}\text{, }z\in \mathbf{R}\text{, }v\in \mathbf{R}.
\end{array}
\right.
\]

\item[(A2)]  There exist some constant $c>0$ and $0<\alpha <1$ such that for
all $\left( t,x,y_i,z_i,v\right) \in \left[ 0,T\right] \times \mathbf{%
R\times R\times R\times R},$ ($i=1,2$),
\[
\left\{
\begin{array}{l}
\left| f\left( t,x,y_1,z_1,v\right) -f\left( t,x,y_2,z_2,v\right) \right|
^2\leq c\left( \left| y_1-y_2\right| ^2+\left| z_1-z_2\right| ^2\right) , \\
\left| g\left( t,x,y_1,z_1\right) -g\left( t,x,y_2,z_2\right) \right| ^2\leq
c\left| y_1-y_2\right| ^2+\alpha \left| z_1-z_2\right| ^2.
\end{array}
\right.
\]
\end{enumerate}

Let $\mathcal{U}_{ad}$ be an admissible control set. The optimal control
problem of SPDE (5.1) is to find an optimal control , such that
\[
J\left( v_{\left( \cdot \right) }^{*}\right) \doteq \inf\limits_{v_{\left(
\cdot \right) }\in \mathcal{U}_{ad}}J\left( v_{\left( \cdot \right) }\right)
,
\]
where $J\left( v_{\left( \cdot \right) }\right) $ is its cost function as
follows:
\begin{equation}
J\left( v_{\left( \cdot \right) }\right) =\mathbf{E}\left[ \int_0^Tl\left(
s,x,u\left( s,x\right) ,\left( \nabla u\sigma \right) \left( s,x\right)
,v_s\right) \text{d}s+\gamma \left( u\left( 0,x\right) \right) \right] .
\tag{5.2}
\end{equation}
Here we assume $l$ and $\gamma $ satisfy (H4). We can transform the optimal
control problem of SPDE (5.1) into one of the following FBDSDE with control
variable:
\begin{equation}
\left\{
\begin{array}{l}
X_s^{t,x}=x+\int_t^sb\left( X_r^{t,x},v_r\right) \text{d}r+\int_t^s\sigma
\left( X_r^{t,x}\right) \text{d}W_r, \\
Y_s^{t,x}=\tilde h\left( X_T^{t,x}\right) +\int_s^Tf\left(
r,X_r^{t,x},Y_r^{t,x},Z_r^{t,x},v_r\right) \text{d}r+\int_s^Tg\left(
r,X_r^{t,x},Y_r^{t,x},Z_r^{t,x}\right) \text{\^d}B_r \\
\quad \quad -\int_s^TZ_r^{t,x}\text{d}W_r,\quad \quad 0\leq t\leq s\leq T,
\end{array}
\right.  \tag{5.3}
\end{equation}
where $\left( X_{\left( \cdot \right) }^{t,x},Y_{\left( \cdot \right)
}^{t,x},Z_{\left( \cdot \right) }^{t,x},v_{\left( \cdot \right) }\right) \in
\mathbf{R\times R\times R\times R}$, $x\in \mathbf{R}$. The corresponding
optimal control problem of FBDSDE (5.3) is to find an optimal control $%
v_{\left( \cdot \right) }^{*}\in \mathcal{U}_{ad}$, such that
\[
J\left( v_{\left( \cdot \right) }^{*}\right) \doteq \inf\limits_{v_{\left(
\cdot \right) }\in \mathcal{U}_{ad}}J\left( v_{\left( \cdot \right) }\right)
,
\]
where $J\left( v_{\left( \cdot \right) }\right) $ is the cost function same
as (5.2):
\[
J\left( v_{\left( \cdot \right) }^{*}\right) \doteq \inf\limits_{v_{\left(
\cdot \right) }\in \mathcal{U}_{ad}}J\left( v_{\left( \cdot \right) }\right)
,
\]

Now we consider the following adjoint FBDSDEs involving the four unknown
processes $\left( p_t,q_t,k_t,h_t\right) $:
\begin{equation}
\left\{
\begin{array}{l}
\text{d}p_t=\left( f_Yp_t+g_Yk_t-l_Y\right) \text{d}t+\left(
f_Zp_t-g_Zk_t-l_Z\right) \text{d}W_t-k_t\text{\^d}B_t, \\
\text{d}q_t=\left( f_Xp_t-b_Xq_t+g_Xk_t-\sigma _Xh_t-l_X\right) \text{d}t+h_t%
\text{d}W_t, \\
p_0=-\gamma _Y\left( Y_0\right) ,\quad q_T=-\tilde h_X\left( X_T\right) p_T,%
\text{\qquad }0\leq t\leq T.
\end{array}
\right.  \tag{5.4}
\end{equation}
It is easy to see that the first equation of (5.4) is a ``forward'' BDSDE,
so it is uniquely solvable by virtue of the result in [14]. The second
equation of (5.4) is a standard BSDE, so it is uniquely solvable by virtue
of the result in [13]. Therefore we know that (5.4) has a unique solution $%
\left( p_{\left( \cdot \right) },q_{\left( \cdot \right) },k_{\left( \cdot
\right) },h_{\left( \cdot \right) }\right) \in M^2\left( 0,T;\mathbf{R\times
R\times R\times R}\right) $. Define the Hamilton function as follows:
\begin{eqnarray*}
\bar H\left( t,X,Y,Z,v,p,q,k,h\right) &=&H\left( t,X,Y,0,Z,v,p,q,k,h\right)
\\
\ &=&l\left( t,X,Y,Z,v\right) -k\cdot g\left( t,X,Y,Z\right) \\
&&\ +q\cdot b\left( X,v\right) -p\cdot f\left( t,X,Y,Z,v\right) +h\cdot
\sigma \left( X\right) .
\end{eqnarray*}
\begin{equation}
\tag{5.5}
\end{equation}
We now formulate a maximum principle for the optimal control system of (5.3).

\begin{theorem}
Suppose (A1)-(A2) hold. Let $\left( X_{\left( \cdot \right) },Y_{\left(
\cdot \right) },Z_{\left( \cdot \right) },u_{\left( \cdot \right) }\right) $
be an optimal control and its corresponding trajectory of (5.3), $\left(
p_{\left( \cdot \right) },q_{\left( \cdot \right) },k_{\left( \cdot \right)
},h_{\left( \cdot \right) }\right) $ be the solution of (5.4). Then the
maximum principle holds, that is, for $t\in \left[ 0,T\right] $, $\forall
v\in \mathcal{U},$
\[
\bar H\left( t,X_t,Y_t,Z_t,v,p_t,q_t,k_t,h_t\right) \geq \bar H\left(
t,X_t,Y_t,Z_t,v_t^{*},p_t,q_t,k_t,h_t\right) ,\text{ a.e., a.s..}
\]
\end{theorem}

\noindent \textbf{Proof.}\quad\noindent Noting that the forward equation of
(5.3) is independent of the backward one, we easily know that it is uniquely
solvable. It is straightforward to use the same arguments in Section 3 to
obtain the desired results. We omit the detailed proof.\quad $\Box$

From the results in [14], we easily have the following propositions.

\begin{proposition}
For any given admissible control $v_{\left( \cdot \right) },$ we assume (A1)
and (A2) hold. Then (5.3) has a unique solution $\left( X_{\left( \cdot
\right) }^{t,x},Y_{\left( \cdot \right) }^{t,x},Z_{\left( \cdot \right)
}^{t,x}\right) \in M^2\left( 0,T;\mathbf{R\times R\times R}\right) $.
\end{proposition}

\begin{proposition}
For any given admissible control $v_{\left( \cdot \right) },$ we assume (A1)
and (A2) hold. Let $\left\{ u\left( t,x\right) ;0\leq t\leq T,x\in \mathbf{R}%
\right\} $ be a random field such that $u\left( t,x\right) $ is $\mathcal{F}%
_{t,T}^B$-measurable for each $\left( t,x\right) ,$ $u\in C^{0,2}\left(
\left[ 0,T\right] \times \mathbf{R};\mathbf{R}\right) $ a.s., and $u$
satisfies SPDE (5.1). Then $u\left( t,x\right) =Y_t^{t,x}.$
\end{proposition}

\begin{proposition}
For any given admissible control $v_{\left( \cdot \right) },$ we assume (A1)
and (A2) hold. Then $\left\{ u\left( t,x\right) =Y_t^{t,x};0\leq t\leq
T,x\in \mathbf{R}\right\} $ is a unique classical solution of SPDE (5.1).
\end{proposition}

Set the Hamilton function
\begin{eqnarray*}
\bar H\left( t,x,u,\nabla u\sigma ,v,p,q,k,h\right) &=&l\left( t,x,u,\nabla
u\sigma ,v\right) -k\cdot g\left( t,x,u,\nabla u\sigma \right) \\
&&\ \ +q\cdot b\left( x,v\right) -p\cdot f\left( t,x,u,\nabla u\sigma
,v\right) +h\cdot \sigma \left( x\right) .
\end{eqnarray*}
Now we can state the maximum principle for the optimal control problem of
SPDE (5.1).

\begin{theorem}
Suppose $u\left( t,x\right) $ is the optimal solution of SPDE (5.1)
corresponding to the optimal control $v_{\left( \cdot \right) }^{*}$ of
(5.1). Then we have, for any $v\in \mathcal{U}$ and $t\in \left[ 0,T\right] ,
$ $x\in \mathbf{R,}$
\begin{eqnarray*}
&&\ \ \bar H\left( t,x,u\left( t,x\right) ,\left( \nabla u\sigma \right)
\left( t,x\right) ,v,p_t,q_t,k_t,h_t\right) \\
\ &\geq &\bar H\left( t,x,u\left( t,x\right) ,\left( \nabla u\sigma \right)
\left( t,x\right) ,v_t^{*},p_t,q_t,k_t,h_t\right) ,\text{ a.e., a.s.}
\end{eqnarray*}
\end{theorem}

\noindent\textbf{Proof.}\quad By virtue of Proposition 13, 14 and 15, the
optimal control problem of SPDE (5.1) can be transformed into the one of
FBDSDE (5.3). Hence, from Theorem 12, the desired result is easily obtained.
\quad $\Box$

\textbf{Remark.}\quad In Section 5, we study the optimal control problem of
a kind of quasilinear SPDE which was similar to the SPDE considered by \O
ksendal in [12]. It is worth mentioning that the quasilinear SPDEs in [12]
can also be related to a class of partially coupled FBDSDEs. Consequently
the results in [12] can be obtained by the approach of FBDSDEs.

\section{Linear quadratic nonzero sum doubly stochastic differential games}

In this section, we investigate linear quadratic non zero sum doubly
stochastic differential games problem. Under the framework of
uniqueness and existence result introduced above, we improve similar
result in Hamadene [11] and Wu [27]. For natational simplification,
we only consider two players, which is similar for $n$ players. Now
the control system is
\begin{equation}
\left\{
\begin{array}{l}
\text{d}x_t^v=\left[ Ax_t^v+B^1v_t^1+B_t^2v_t^2+Ck_t^v+\alpha
_t\right]
\text{d}t+\left[ Dx_t^v+Ek_t^v+\beta _t\right] \text{d}W_t-k_t^v\text{\^d}%
B_t, \\
x_0^v=a,\qquad t\in \left[ 0,T\right] ,
\end{array}
\right.  \tag{6.1}
\end{equation}

where $A,$ $C,$ $D$ and $E$ are $n\times n$ bounded matrices,
further, $E$ satisfies $0<\left| E\right| <1,$ $v_t^1$ and $v_t^2$,
$t\in \left[
0,T\right] ,$ are two admissible control processes, that is ${\cal F}_t$%
-progressively measurable square integrable processes taking values in $R^k$%
. $B^1$ and $B^2$ are $n\times k$ bounded matrices. $\alpha _t$ and
$\beta _t $ are two adapted squares-integrable processes. We denote
by
\begin{equation}
\left\{
\begin{array}{c}
J^1\left( v\left( \cdot \right) \right) =\frac 12{\bf E}\left[
\int_0^T\left( \left\langle R^1x_t^v,x_t^v\right\rangle
+\left\langle N^1v_t^1,v_t^1\right\rangle +\left\langle
P^1k_t^v,k_t^v\right\rangle
\right) \text{d}t+\left\langle Q^1x_T^v,x_T^v\right\rangle \right] , \\
J^2\left( v\left( \cdot \right) \right) =\frac 12{\bf E}\left[
\int_0^T\left( \left\langle R^2x_t^v,x_t^v\right\rangle
+\left\langle N^2v_t^2,v_t^2\right\rangle +\left\langle
P^2k_t^v,k_t^v\right\rangle \right) \text{d}t+\left\langle
Q^2x_T^v,x_T^v\right\rangle \right] .
\end{array}
\right.  \tag{6.2}
\end{equation}

Here $Q^i,$ $R^i,$ and $P^i,$ $i=1,2$, are $n\times n$ nonnegative
symmetric bounded matrices, $N^1$ and $N^2$ are $k\times k$ positive
symmetric bounded matrices and inverses $\left( N^1\right) ^{-1},$
$\left( N^2\right) ^{-1}$ are also bounded. We denote $v\left( \cdot
\right) =\left( v^1\left( \cdot \right) ,v^2\left( \cdot \right)
\right) .$ The problem is to find $\left( u^1\left( \cdot \right)
,u^2\left( \cdot \right) \right) \in R^k\times R^k$ which is called
Nash equilibrium point for the game, such that
\begin{equation}
\left\{
\begin{array}{c}
J^1\left( u^1\left( \cdot \right) ,u^2\left( \cdot \right) \right)
\leq J^1\left( v^1\left( \cdot \right) ,u^2\left( \cdot \right)
\right) ,\quad
v^1\left( \cdot \right) \in R^k; \\
J^2\left( u^1\left( \cdot \right) ,u^2\left( \cdot \right) \right)
\leq J^2\left( u^1\left( \cdot \right) ,v^2\left( \cdot \right)
\right) ,\quad v^1\left( \cdot \right) \in R^k.
\end{array}
\right.  \tag{6.3}
\end{equation}

Note that the actions of the two players are described by a
classical BDSDE in which we indicates that the players should make
some strategy to overcome the disturbed information. In order to
introduce the main result, we need the following assumptions:
\begin{equation}
\left\{
\begin{array}{c}
B^i\left( N^i\right) ^{-1}\left( B^i\right) ^TA^T=A^TB^i\left(
N^i\right)
^{-1}\left( B^i\right) ^T,\quad i=1,2, \\
B^i\left( N^i\right) ^{-1}\left( B^i\right) ^TC^T=C^TB^i\left(
N^i\right)
^{-1}\left( B^i\right) ^T,\quad i=1,2, \\
B^i\left( N^i\right) ^{-1}\left( B^i\right) ^TD^T=D^TB^i\left(
N^i\right)
^{-1}\left( B^i\right) ^T,\quad i=1,2, \\
B^i\left( N^i\right) ^{-1}\left( B^i\right) ^TE^T=E^TB^i\left(
N^i\right)
^{-1}\left( B^i\right) ^T,\quad i=1,2, \\
B^i\left( N^i\right) ^{-1}\left( B^i\right) ^TP^1=P^1B^i\left(
N^i\right)
^{-1}\left( B^i\right) ^T,\quad i=1,2, \\
B^i\left( N^i\right) ^{-1}\left( B^i\right) ^TP^2=P^2B^i\left(
N^i\right) ^{-1}\left( B^i\right) ^T,\quad i=1,2,
\end{array}
\right.  \tag{6.4}
\end{equation}

Next we give an explicit form of Nash equilibrium point by virtue of
solutions of linear FBDSDEs. Hence we have a following theorem.

\begin{theorem}
The pair of function
\[
\left\{
\begin{array}{l}
u_t^1=-\left( N^1\right) ^{-1}\left( B^1\right) ^Ty_t^1, \\
u_t^2=-\left( N^1\right) ^{-1}\left( B^1\right) ^Ty_t^2,\quad t\in
\left[ 0,T\right] ,
\end{array}
\right.
\]

is one Nash equilibrium point for the above game problem, where
$\left( x_t,y_t^1,y_t^2,k_t^1,k_t^2,h_t^1,h_t^2\right) $ is the
solution of the following differential dimensinal FBDSDEs:
\begin{equation}
\left\{
\begin{array}{l}
\text{d}x_t=\left[ Ax_t-B^1\left( N^1\right) ^{-1}\left( B^1\right)
^Ty_t^1-B^2\left( N^2\right) ^{-1}\left( B^2\right) ^Ty_t^2+\alpha
_t\right]
\text{d}t \\
\qquad \quad \left[ Cx_t+\beta _t\right] \text{d}W_t-k_t\text{\^d}B_t, \\
\text{d}y_t^1=-\left[ Ay_t^1+D^Th_t^1+R^1x_t\right] \text{d}t-\left(
C^Ty_t^1+E^Th_t^1+P^1k_t\right) \text{\^d}B_t+h_t^1\text{d}W_t, \\
\text{d}y_t^2=-\left[ Ay_t^2+D^Th_t^2+R^2x_t\right] \text{d}t-\left(
C^Ty_t^2+E^Th_t^2+P^2k_t\right) \text{\^d}B_t+h_t^2\text{d}W_t, \\
x_0=a,\quad y_T^1=Q^1x_T,\quad y_T^2=Q^2x_T.
\end{array}
\right.  \tag{6.5}
\end{equation}
\end{theorem}

\noindent \textbf{Proof of Lemma 2.}\quad At the beginning, we prove
the existence of the solution of (6.5). Consider the following
FBDSDEs:
\begin{equation}
\left\{
\begin{array}{l}
\text{d}X_t=\left( AX_t-Y_t+\alpha _t\right) \text{d}t+\left[
CX_t+\beta
_t\right] \text{d}W_t-K_t\text{\^d}B_t, \\
\text{d}Y_t=-\left( A^TY_t+\left( \left( B^1\left( N^1\right)
^{-1}\right) \left( B^1\right) ^TR^1+\left( B^1\left( N^1\right)
^{-1}\right) \left(
B^1\right) ^TR^2\right) X_t+D^TH_t\right) \text{d}t \\
\qquad -\left[ C^TY_t+E^TH_t+PK_t\right] \text{\^d}B_t+H_t\text{d}W_t, \\
X_0=a,\quad Y_T=\left[ \left( B^1\left( N^1\right) ^{-1}\right)
\left( B^1\right) ^TR^1+\left( B^1\left( N^1\right) ^{-1}\right)
\left( B^1\right) ^TR^1\right] X_T.
\end{array}
\right.  \tag{6.6}
\end{equation}

Apparently, if the $\left(
x_t,y_t^1,y_t^2,k_t^1,k_t^2,h_t^1,h_t^2\right) $ is the solution of
(6.5), then $\left( X_t,Z_t,Y_t\right) $ satisfies the FBDSDEs (6.6)
with (6.4). Here
\[
\left\{
\begin{array}{l}
X_t=x_t, \\
K_t=k_t, \\
Y_t=B^1\left( N^1\right) ^{-1}\left( B^1\right) ^Ty_t^1+B^2\left(
N^2\right)
^{-1}\left( B^2\right) ^Ty_t^2, \\
H_t=B^1\left( N^1\right) ^{-1}\left( B^1\right) ^Th_t^1+B^2\left(
N^2\right)
^{-1}\left( B^2\right) ^Th_t^2, \\
P=P^1B^1\left( N^1\right) ^{-1}\left( B^1\right) ^T+P^2B^2\left(
N^2\right) ^{-1}\left( B^2\right) ^T.
\end{array}
\right.
\]

As matter of fact, it is easy to check that there exits a unique solution $%
\left( X_t,Y_t,K_t,H_t\right) $ of (6.6) according to Proposition 2.
Hence we can first solve the FBDSDEs (6.6) to get solution $\left(
X_t,K_t\right) $ which, obviously, is the forward solution $\left(
x_t,k_t\right) $ of (6.5), then $\left( y_t^1,h_t^1\right) $ and
$\left( y_t^2,h_t^2\right) $ are obtained. Now consider the
following classical backward doubly stochastic differential
equations (BDSDEs in short) with four unknown processes $\left(
y_t^1,y_t^2,h_t^1,h_t^2\right) $:
\[
\left\{
\begin{array}{l}
\text{d}y_t^1=-\left[ A^Ty_t^1+C^Th_t^1+R^1X_t\right]
\text{d}t-\left(
C^Ty_t^1+E^Th_t^1+P^1k_t\right) \text{\^d}B_t+h_t^1\text{d}W_t, \\
\text{d}y_t^2=-\left[ A^Ty_t^2+C^Th_t^2+R^2X_t\right]
\text{d}t-\left(
C^Ty_t^2+E^Th_t^2+P^2k_t\right) \text{\^d}B_t+h_t^2\text{d}W_t, \\
y_t^1=Q^1X_T,\qquad y_t^2=Q^2X_T.
\end{array}
\right.
\]

Set
\[
\left\{
\begin{array}{c}
\hat Y_t=B^1\left( N^1\right) ^{-1}\left( B^1\right)
^Ty_t^1+B^2\left(
N^2\right) ^{-1}\left( B^2\right) ^Ty_t^2, \\
\hat H_t=B^1\left( N^1\right) ^{-1}\left( B^1\right)
^Th_t^1+B^2\left( N^2\right) ^{-1}\left( B^2\right) ^Th_t^2,
\end{array}
\right.
\]

after simple computation we have
\[
\left\{
\begin{array}{l}
\text{d}\hat Y_t=-\left[ A^T\hat Y_t+\left( B^1\left( N^1\right)
^{-1}\left( B^1\right) ^TR^1+B^2\left( N^2\right) ^{-1}\left(
B^2\right) ^TR^2\right)
X_t+D^T\hat H_t\right] \text{d}t \\
\qquad -\left[ C^TY_t+E^TH_t+PK_t\right] \text{\^d}B_t+\hat
H_t\text{d}W_t,
\\
\hat Y_T=\left[ B^1\left( N^1\right) ^{-1}\left( B^1\right)
^TR^1+B^2\left( N^2\right) ^{-1}\left( B^2\right) ^TR^2\right] X_T.
\end{array}
\right.
\]

Now fixing $\left\{ X_t\right\} _{t\geq 0}$, and thanks to $0<\left|
E\right| <1,$ due to the existence and uniqueness of solution of
BDSDE, we immediately have
\[
\left\{
\begin{array}{c}
Y_t=\hat Y_t=B^1\left( N^1\right) ^{-1}\left( B^1\right)
^Ty_t^1+B^2\left(
N^2\right) ^{-1}\left( B^2\right) ^Ty_t^2, \\
H_t=\hat H_t=B^1\left( N^1\right) ^{-1}\left( B^1\right)
^Th_t^1+B^2\left( N^2\right) ^{-1}\left( B^2\right) ^Th_t^2,
\end{array}
\right.
\]

No doubt, $\left( X_t,Y_t,K_t,H_t\right) $ satisfies the FBDSDEs
(6.6) and is the unique solution. Therefore $\left(
x_t,y_t^1,y_t^2,k_t^1,k_t^2,h_t^1,h_t^2\right) $ is the solution of
FBDSDEs (6.5). From now on we prove $\left( u^1\left( t\right)
,u^2\left( t\right) \right) $ is one Nash equilibrium point for our
nonzero sum game problem. For that it suffices that
\[
J^1\left( u^1\left( \cdot \right) ,u^2\left( \cdot \right) \right)
\leq J^1\left( v^1\left( \cdot \right) ,u^2\left( \cdot \right)
\right) ,\qquad \forall v^1\left( \cdot \right) \in R^k.
\]

It is similar to give the other inequality by the same argument.
Next we give the control system by $x_t^{v^1}$:
\[
\left\{
\begin{array}{l}
\text{d}x_t^{v^1}=\left[
Ax_t^{v^1}+B^1v_t^1+B^2u_t^2+Ck_t^{v^1}+\alpha
_t\right] \text{d}t+\left[ Cx_t^{v^1}+\beta _t\right] \text{d}W_t-k_t^{v^1}%
\text{d}B_t, \\
x_0=a,\qquad t\in \left[ 0,T\right] ,
\end{array}
\right.
\]

\begin{eqnarray*}
&&\ J^1\left( v^1\left( \cdot \right) ,u^2\left( \cdot \right)
\right)
-J^1\left( u^1\left( \cdot \right) ,u^2\left( \cdot \right) \right) \\
&=&\frac 12{\bf E[}\int_0^T\left( \left\langle
R^1x_t^{v^1},x_t^{v^1}\right\rangle -\left\langle
R^1x_t,x_t\right\rangle
+\left\langle N^1v_t^1,v_t^1\right\rangle \right. \\
&&\ \qquad \left. -\left\langle N^1u_t^1,u_t^1\right\rangle
+\left\langle P^1k_t^{v^1},k_t^{v^1}\right\rangle -\left\langle
P^1k_t,k_t\right\rangle
\right) \text{d}t \\
&&\ \ \ \qquad +\left\langle Q^1x_T^{v^1},x_T^{v^1}\right\rangle
-\left\langle Q^1x_T,x_T\right\rangle ] \\
&=&\frac 12{\bf E[}\int_0^T\left( \left\langle R^1\left(
x_t^{v^1}-x_t\right) ,x_t^{v^1}-x_t\right\rangle \right. \\
&&\qquad +\left\langle N^1\left( v_t^1-u_t^1\right)
,v_t^1-u_t^1\right\rangle
\\
&&\qquad +\left\langle P^1\left( k_t^{v^1}-k_t\right)
,k_t^{v^1}-k_t\right\rangle +2\left\langle
R^1x_t,x_t^{v^1}-x_t\right\rangle
\\
&&\qquad \left. +2\left\langle N^1u_t^1,v_t^1-u_t^1\right\rangle
+2\left\langle P^1k_t,k_t^{v^1}-k_t\right\rangle \right) \text{d}t \\
&&\qquad +\left\langle Q^1\left( x_T^{v^1}-x_T\right)
,x_T^{v^1}-x_T\right\rangle \\
&&\qquad +2\left( Q^1x_T,x_T^{v^1}-x_T\right) ].
\end{eqnarray*}

Note that
\[
Q^1x_T=y_T^1.
\]

We apply It\^o's formula to $\left\langle
x_T^{v^1}-x_T,y_T^1\right\rangle $ on the $\left[ 0,T\right] $ and
get
\begin{eqnarray*}
{\bf E}\left\langle x_T^{v^1}-x_T,y_T^1\right\rangle &=&{\bf E}%
\int_0^T\left( -\left\langle R^1x_t,\left( x_t^{v^1}-x_t\right)
\right\rangle +\left\langle B^1\left( v_t^1-u_t^1\right)
,y_t^1\right\rangle
\right. \\
&&\left. -\left\langle P^1k_t,k_t^{v^1}-k_t\right\rangle \right)
\text{d}t.
\end{eqnarray*}

Under the assumption $R^1$, $Q^1$ and $P^1$ being nonnegative, $N^1$
being positive, and symmetry of $B^1,$ we have
\begin{eqnarray*}
&&J^1\left( v^1\left( \cdot \right) ,u^2\left( \cdot \right) \right)
-J^1\left( u^1\left( \cdot \right) ,u^2\left( \cdot \right) \right) \\
&\geq &{\bf E}\int_0^T\left( \left\langle
N^1u_t^1,v_t^1-u_t^1\right\rangle
+\left\langle B^1\left( v_t^1-u_t^1\right) ,y_t^1\right\rangle \right) \text{%
d}t \\
&=&{\bf E}\int_0^T\left( \left\langle -N^1\left( N^1\right)
^{-1}\left( B^1\right) ^Ty_t^1,v_t^1-u_t^1\right\rangle
+\left\langle \left( B^1\right)
^Ty_t^1,v_t^1-u_t^1\right\rangle \right) \text{d}t \\
&=&0.
\end{eqnarray*}

Lastly, we claim that
\[
\left\{
\begin{array}{l}
u_t^1=-\left( N^1\right) ^{-1}\left( B^1\right) ^Ty_t^1, \\
u_t^2=-\left( N^1\right) ^{-1}\left( B^1\right) ^Ty_t^1,\quad t\in
\left[ 0,T\right] ,
\end{array}
\right.
\]

that is, $\left( u_t^1,u_t^2\right) $ is one Nash equilibrium point
for our nonzero sum doubly stochastic game problem.\quad $\Box$

\textbf{Remark 2}\ As matter of fact, in Theorem 17, we use the
adjoint equation, the idea is the same as in Theorem 10. Besides,
the results of this section are clear and easy to understand. They
can be applied in practice directly.

{\bf Acknowledgment}: The authors would like to thank the referees
for their helpful comments and suggestions.

\end{document}